\documentclass[11pt]{article}
\usepackage[margin=1in]{geometry}
\usepackage{marvosym}   % \Letter, the corresponding-author mark
\newenvironment{keywords}{\par\medskip\noindent\textbf{Key words.} }{\par}
\newenvironment{AMS}{\par\medskip\noindent\textbf{AMS subject classifications.} }{\par\medskip}
\AtBeginDocument{\numberwithin{equation}{section}}   % per-section equation numbers, as in the SIAM class
\usepackage{amsmath}
\usepackage{amsfonts}
\usepackage{amssymb}
\usepackage{graphicx}
\DeclareGraphicsExtensions{.pdf,.png,.jpg}
\usepackage{bm}
\usepackage{xcolor} 
\usepackage{graphicx,epstopdf,tikz}
 \DeclareMathAlphabet{\itbf}{OML}{cmm}{b}{it}
\DeclareMathAlphabet\mathbfcal{OMS}{cmsy}{b}{n}

\renewcommand{\hat}{\widehat}
\renewcommand{\tilde}{\widetilde}

\newcommand{\RR}{\mathbb{R}}

\def\om{\omega}
\def\la{\lambda}

\def\bx{{{\itbf x}}}

\def\bz{{{\itbf z}}}
\def\bu{{{\itbf u}}}

\def\btheta{\boldsymbol{\theta}}
\def\bg{{\itbf g}}

\def\bU{{\itbf U}}
\def\bJ{{\itbf J}}

\def\bv{{\bf v}}
\def\bet{{\boldsymbol{\eta}}}
\def\ss{{(s)}}

\def\bV{{\itbf V}}
\def\bR{{\itbf R}}

\def\bS{\boldsymbol{\mathbb{S}}}

\def\bM{\boldsymbol{\mathbb{M}}}

\def\bA{\mathbfcal{A}}  % was \boldsymbol{\cal A}; switched to predeclared \mathbfcal to avoid bm's math-alphabet overflow

\newcommand{\cG}{\mathcal{G}}

\def\bGa{{\boldsymbol{\Gamma}}}
\def\bcM{{\boldsymbol{\mathcal M}}}
\def\bcD{{\boldsymbol{\mathcal D}}}

\def\bPhi{{\boldsymbol{\Phi}}}

\def\FWI{{\scalebox{0.5}[0.4]{FWI}}}

\def\RM{{\scalebox{0.5}[0.4]{ROM}}}
\def\Gal{{\scalebox{0.5}[0.4]{GAL}}}

\def\ML{{\scalebox{0.5}[0.4]{ML}}}

\def\cF{\mathcal{F}}
\def\cT{\mathcal{T}}
\def\12{{\frac{1}{2}}}

\usepackage{pgfplots}
\usepackage{tikz}
\usepackage{algorithm}
\usepackage[noend]{algpseudocode}
\usepackage{booktabs}
\usepackage{multirow}
\usepackage{array}

\newtheorem{prop}{Proposition}[section]

\usepackage[normalem]{ulem}

\begin{document}
\renewcommand{\thefootnote}{\fnsymbol{footnote}}   % author footnotes use symbols, as in the SIAM class

\title{ROMNet: a hybrid reduced order modeling and machine learning approach  to waveform inversion}
\date{}   % no date, as in the SIAM version
 \author{Liliana Borcea\footnotemark[1] \quad Alexander Mamonov\footnotemark[2] \quad 
Kui Ren\footnotemark[1] \quad Haizhao Yang\footnotemark[5] \quad  Chugang Yi\footnotemark[5]\mbox{\hspace{7pt}}{\scriptsize\Letter{}} }

\maketitle

 %%%%%%%%%%%%%%%%%%%%%%%%%%%%%
 % footnotes for the addresses

{\renewcommand{\thefootnote}{}\footnotetext{Authors are listed in alphabetical order.}}
\footnotetext[1]{Applied Physics and Applied Mathematics, Columbia University, New York, NY,  10027. {\tt lb3539@columbia.edu} and 
{\tt kr2002@columbia.edu}}
\footnotetext[2]{Department of Mathematics, University of Houston,  Houston, TX 77204-3008. {\tt mamonov@math.uh.edu}}
\footnotetext[5]{Department of Mathematics, University of Maryland, College Park, MD 20742. {\tt chugang@umd.edu} and 
{\tt hzyang@umd.edu}}
{\renewcommand{\thefootnote}{{\scriptsize\Letter}}\footnotetext[1]{Corresponding author ({\tt chugang@umd.edu}).}}
\renewcommand{\thefootnote}{\arabic{footnote}}\setcounter{footnote}{0}

\begin{abstract} 
Waveform inversion seeks to estimate the wave speed of a heterogeneous, inaccessible medium, from time-resolved measurements of the waves at user controlled sensors. We consider this inverse problem for acoustic waves and an active array of source/receiver sensors that emit probing signals and measure the generated pressure waves. The forward map, from the wave speed to the measurements, is nonlinear and oscillatory. The oscillations cause cycle skipping,  the main impediment to using the standard, nonlinear least-squares data fitting formulation, known as full waveform inversion (FWI). A recently introduced alternative waveform inversion approach computes from the measurements an algebraic surrogate of the wave operator, a reduced order model (ROM) matrix, which is then used to estimate the wave speed. The mapping from the measurements to the ROM is nonlinear, but well understood. It is computed efficiently, in a non-iterative manner. The nonlinear mapping from the ROM to the wave speed is less understood, and its approximation involves time-consuming optimization. Our goal in this paper is to use a neural network to map the ROM matrix to a nearby one, that has a simpler and explicit dependence on the wave speed. This simplifies and reduces the computational cost of the ROM-based waveform inversion. We introduce the methodology, called ROMNet, and test it with numerical simulations, using two training data sets: The first set consists of random media with variations of the wave speed modeled by a superposition of Gaussians with random amplitudes and standard deviations. The second is the publicly available GeoFWI dataset introduced for benchmarking FWI using deep learning. We compare the performance of ROMNet with the direct ROM-based inversion and with two representative deep learning approaches to FWI: ``Fourier-DeepONet" and ``InversionNet".
\end{abstract}

% Note that keywords are not normally used for peer review papers.
\begin{keywords}
waveform inversion, data driven reduced order modeling, machine learning.
\end{keywords}

\begin{AMS}
35-R30, 65-M32,  86-A22
\end{AMS}

\section{Introduction}
Waveform inversion is an inverse problem for the wave equation, which seeks to estimate a heterogeneous, inaccessible medium, from measurements of the wave field gathered at user controlled sensors. It is a widely studied problem with a diverse range of applications in medical ultrasound, airborne or ground based radar, exploration geophysics, nondestructive evaluation of materials, etc.  Depending on the application, the wave equation models sound, elastic or electromagnetic waves. We 
work with  sound waves in a medium with constant density and unknown wave speed $c(\bx)$. The wave field is the acoustic pressure $p^{(s)}(t,\bx)$, labeled by the source index $s$. It satisfies the wave equation
\begin{align}
\left[\partial_t^2 - c^2(\bx) \Delta \right] p^{(s)}(t,\bx) &= f(t) \delta(\bx-\bx_s), \quad t \in \RR, ~~ \bx \in \RR^d, 
\label{eq:I1} 
\end{align}
with  initial condition 
\begin{align}
p^{(s)}(t,\bx) &= 0, \quad t < -t_f.
\label{eq:I2}
\end{align}
We state the problem in $\RR^d$, with dimension $d \in \{1,2,3\}$. The methodology applies to any dimension, but all the numerical simulations are  for $d=2$.

The data consist of wave measurements  gathered by  an active array of $m$ sensors located at  $\bx_s$, for $1 \le s \le m$. The array uses one source 
sensor at a time, hence the forcing in equation \eqref{eq:I1}.  It emits a probing pulse $f(t)$ supported in the interval $(-t_f,t_f)$ and then measures the generated pressure over some finite time interval ${\mathcal T}$. 
Because the array has a  one sided view of the medium, we refer to the measurements as ``back-scattering data".
The waveform inversion problem is to estimate $c(\bx)$ from the array response  back-scattering data matrix $\bcM(t)\in \RR^{m \times m}$ with entries
\begin{equation}
\mathcal{M}^{(r,s)}(t) = p^\ss(t,\bx_r), \qquad 1\le r,s\le m, \quad t \in \mathcal{T}. \label{eq:I3}
\end{equation}

There is extensive literature on waveform inversion, much of it developed in tandem with applications. The methodology is divided roughly in two categories: qualitative, imaging methods and quantitative methods. Qualitative methods address a simpler problem: Localize the reflectivity of the medium, modeled by jump discontinuities of $c(\bx)$, in a known, smooth, often homogeneous background \cite{biondi20063d,symes2008migration,jakowatz2012spotlight,cheney2009fundamentals,cakoni2014qualitative}. 
Quantitative methods seek to determine the wave speed function i.e., both the smooth and rough parts of $c(\bx)$ \cite{tarantola2005inverse,VirieuxOperto2009}. 

The analysis and implementation of both qualitative and quantitative  methods depends on the frequency content of the probing waves. For example, the linear sampling 
method  \cite{colton1996simple,colton1997simple}, the factorization method \cite{kirsch1998characterization,kirsch1999factorization} and the related MUSIC (Multiple Signal Classification) algorithm
\cite{stoica1989music,cheney2001linear} are developed in the frequency domain, for  time harmonic waves. They are qualitative methods studied in detail \cite{cakoni2014qualitative,kirsch2007factorization,ammari2005music,borcea2016robust}, because they are non-iterative and  take multiple scattering into account. MUSIC is generally used for localizing  point-like reflectors, while linear sampling and the factorization method 
work with general targets. However, they require wide angle views of such targets, to compensate for the lack of frequency diversity i.e., they are not optimal for back-scattering data. Broad-band time domain qualitative methods commonly use  single scattering (Born)  approximations with respect to the unknown reflectivity \cite{biondi20063d,symes2008migration,jakowatz2012spotlight,cheney2009fundamentals}. They suffer from artifacts due to neglected multiple scattering effects, but they are robust to high levels of noise, they are flexible with respect to the data acquisition geometry and they can be extended to randomly 
perturbed backgrounds \cite{borcea2006adaptive}. 

Illustrative studies of quantitative methods implemented in the frequency domain are in \cite{bao2015inverse,borges2017high}. They give high resolution 
estimates of $c(\bx)$ from data with very large bandwidth that can be swept systematically from the lowest to the highest frequency, as proposed in \cite{chen1997inverse}.  Such bandwidths are rarely available in applications,  so frequency sweeping is usually not an option.  

Much of the quantitative inversion 
methodology is implemented in the time domain, where the causality of wave propagation is exploited in the computational pursuit of the ``inverse" of the forward map  $c(\bx) \stackrel{\cal F}{\longrightarrow} \bcM(t)$. Note that  $\cF$ is nonlinear and non-invertible in the strict mathematical sense, hence the quotation marks. The problem is generically formulated as a nonlinear least squares data fitting optimization 
\begin{equation}
\mbox{arg}\min_{\hat c \in \mathcal{C}}\,  \mathcal{O}^{\rm FWI}(\, \hat c \, ) + \mbox{regularization}, \qquad  \mathcal{O}^{\rm FWI}(\,\hat c\,) = \int_{\mathcal{T}} dt \, \|{\cal F}[\, \hat c\, ](t) - {\bcM}(t) \|_F^2,
\label{eq:I4}
\end{equation}
where the search speed $\hat c(\bx)$ is in some user defined  space $\mathcal{C}$ and $\| \cdot \|_F^2$ is the matrix Frobenius norm. This formulation is known as 
FWI (Full Waveform Inversion) \cite{tarantola2005inverse,VirieuxOperto2009}. It is a computationally demanding PDE constrained optimization that is flexible with respect to the data acquisition geometry and robust to additive, Gaussian noise. Typically, $\mathcal{C}$ is a high dimensional search space, so FWI is implemented with local optimization methods like Gauss-Newton or gradient descent \cite{gill2019practical}.  These are likely to fail for high frequency, band limited back-scattering data, even when starting from reasonable initial guesses of $c(\bx)$ \cite{VirieuxOperto2009,symes2020wavefield}. The culprit is the ``cycle skipping" phenomenon, that causes many spurious local minima  of $\mathcal{O}^{\rm FWI}(\, \hat c \,)$. Reliable wave speed estimation via local optimization of \eqref{eq:I4} 
requires accurate prior knowledge of the kinematics of the medium, determined by the smooth part of $c(\bx)$, so that the error of the arrival time of the 
predicted waves in $\mathcal{F}[\, \hat c\, ](t)$ is within half a cycle of the oscillation of the data \cite{VirieuxOperto2009,symes2020wavefield}.  This is why the basin of attraction of the global minimum of  $\mathcal{O}^{\rm FWI}(\, \hat c \,)$ shrinks significantly at high frequency \cite{barucq2019priori}.

The forward map restricted to search velocities $\hat{c}(\bx)$ with accurate kinematics has a simpler behavior. Often, it is replaced by its Born approximation, which linearizes ${\cal F}[\, \hat c \, ](t)$  with respect to the rough (reflective) part of $\hat{c}(\bx)$. This leads to a linear least squares data fitting problem that  is the foundation of popular qualitative (imaging) methods: Reverse time (Kirchhoff) migration \cite{biondi20063d,symes2008migration}, matched field imaging \cite{baggeroer1993overview} and filtered backprojection \cite{jakowatz2012spotlight,cheney2009fundamentals}. 

The kinematics is largely unknown in applications involving complex heterogeneous media like earth models in geophysics and extended targets in radar, so the research for mitigating cycle skipping remains an important topic. One  idea is to quantify the data misfit in a different metric than  $L^2(\mathcal{T})$. In particular, metrics from optimal transportation \cite{Engquist2016otfwi,yang2018application,Metivier2018graph} have shown promise in enlarging the basin of attraction of global minima of the resulting objective function. Another  approach is to extend the search space, by adding in a systematic way additional degrees of freedom, like treating the probing signal as unknown \cite{vanLeeuwen2013wri,symes2008migration}. To date, the problem remains largely open and theoretical guarantees that any one approach works do not exist, except for  very simple settings.

A newer point of view, advocated  in \cite{druskin2016direct,borcea2020reduced,Borcea2022rom,borcea2023waveform,borcea2024meets}, 
is  that attempting to ``invert" directly the forward map by data fitting  may not be the best formulation of waveform inversion. These studies  introduce an intermediary mapping, from the data, to a reduced order model (ROM) of wave propagation. Then, they estimate the wave speed from the ROM.  There are 
two types of ROMs: The  ``operator ROM" is a matrix surrogate of the operator $-c^2(\bx) \Delta$  \cite{Borcea2022rom,borcea2024meets}, while the 
``propagator ROM" is a matrix surrogate of the wave propagator operator,  which dictates the evolution of the wave on a uniform time grid \cite{druskin2016direct,borcea2020reduced,borcea2023waveform}. Both  ROMs are constructed carefully, to respect the physics of wave propagation and causality. Mathematically, they  are defined as Galerkin projections of the operators on the space spanned by snapshots of the wave field.  These snapshots cannot be computed in waveform inversion, because the wave speed $c(\bx)$ is unknown. Nevertheless, the ROMs are array data driven, meaning that they are obtained directly from the  response matrix $\bcM(t)$   \cite{druskin2016direct,borcea2020reduced,Borcea2022rom,borcea2023waveform,borcea2024meets}. The mapping from $\bcM(t)$ to the ROMs is nonlinear, but it is well understood and computable in a tractable, non-iterative way, using tools from numerical linear algebra. It is the nonlinear mapping from the ROMs to the unknown $c(\bx)$ that is less understood and requires further study. 

So far, there are two approaches to estimate $c(\bx)$ from the ROMs: The first approach uses the propagator ROM to approximate the wave field at inaccessible points inside the medium. This approximation is consistent with the array data, but not with the solution of the wave equation, unless the search speed equals the true one. Thus, the estimation of $c(\bx)$ is formulated as an iterative optimization that minimizes the discrepancy between the approximated internal wave and the solution of the wave equation \cite{borcea2023waveform,borcea2024meets}. The second approach estimates $c(\bx)$ via iterative minimization of the misfit of the operator ROM. In this paper, we are interested in improving this second approach by training a neural network to map the ROM to a nearby matrix that has an explicit and simple 
dependence on $c(\bx)$. Once this is accomplished, the estimation of $c(\bx)$ becomes easier and less computationally involved.

There is a large and fast-growing literature on machine learning (ML) for waveform inversion~\cite{li2026dlm,saad2024siamesefwi,he2021reparameterized,Ding2022coupling,yan2026deep,YangMa2019}. One approach is to train networks to learn directly the approximate ``inverse" of the forward map $\cF$, like in ``InversionNet" \cite{Wu2019inversionnet} and in ``SeisInvNet" \cite{Li2020seisinvnet}. The training is done on paired synthetic data and wave speed models, and the methods vary in the network designs and choices of grouping data for processing.  Secondary networks have also been added to the training, as in  ``VelocityGAN" \cite{ZhangLin2020velocitygan},  to enhance the resemblance of the predicted wave speed to those in the training set. Extensions to large, three-dimensional problems are considered by ``InversionNet3D" \cite{zeng2021inversionnet3d}. ``Fourier-DeepONet" \cite{zhu2023fourier} is an operator network that adds to the training data the waveforms and variable source locations and frequencies, in order to generalize to source parameters that are not in the training set. Hybrid enhancements of the training process that combine data and physics-guided misfits are proposed in \cite{sun2021physics} and a physics-consistent augmentation of the training data is  in \cite{rojas2020physics}. 

Other ML approaches separate the estimation of the kinematics from that of the reflectivity. One way is to predict the missing low-frequency part of the data~\cite{ovcharenko2019deep,sun2020extrapolated} that facilitates the FWI estimation of  the kinematics. Then, the training can be carried out on the frequency range of the measurements to predict the reflectivity. There are also attempts to reconstruct the smooth and rough parts of $c(\bx)$ separately, {either in a multiscale fashion, from the low to the high frequency content of the data \cite{feng2021multiscale}, or with a scheme driven by the adjoint of the forward map \cite{zhang2020adjoint}}. 

There are various other ML methods, too many to mention here. They all  produce an estimate of $c(\bx)$ at low computational cost, once the training is complete. The risk is that this estimate may be  wrong when the geology, the source/receiver layout and the frequency content are different from those in the training set. This likely  stems from 
the complicated behavior of the $c(\bx)$ to data  map $\cF$. 

In this paper we propose to use ML in a different way, by leveraging the results in~\cite{Borcea2022rom,borcea2024meets} that show that the mapping from $c(\bx)$ to the operator ROM has a better behavior than $\cF$. Explicitly, we introduce the ROMNet approach that simplifies further this mapping, by training a neural network to transform the input operator ROM to a nearby matrix that has an explicit, quadratic dependence on $c(\bx)$. 

The paper is organized as follows. We begin in section \ref{sect:1} with a brief review from \cite{Borcea2022rom,borcea2024meets} of the data-driven computation of the operator ROM and the definition of the mapping $c(\bx) \stackrel{\cF^{\RM}}{\longrightarrow} \bA^{\RM}$. Section \ref{sect:2} introduces and motivates ROMNet. We assess its performance using numerical simulations in Section~\ref{sect:3}. We consider two different training sets for the learning: The first consists of random models of $c(\bx)$ given by a superposition 
of Gaussians with random amplitudes and standard deviations. The second is the GeoFWI data set proposed in \cite{li2026geofwi} as a complement 
to the OpenFWI benchmark~\cite{Deng2022openfwi}. It provides velocity models generated from geological sedimentary and structural rules, which include faults, 
folds, stacked layers, sharp interfaces, and vertical salt columns. The ROMNet results presented in section \ref{sect:3} are for in and out of the training distribution. We compare them with those obtained by the direct ROM approach in \cite{Borcea2022rom,borcea2024meets} and two popular ML approaches: ``Fourier-DeepONet" and ``InversionNet".  We also describe the computational cost savings brought by ROMNet. Some details of the numerical implementation are in appendix \ref{ap:A}. We end in section \ref{sect:4} with a brief summary.

\section{Data-driven operator ROM and its dependence on the wave speed}
\label{sect:1}

Here, we review from \cite{Borcea2022rom,borcea2024meets} the definition and computation of the operator ROM. This is the input of the ROMNet algorithm introduced in the next section. 

The operator ROM is  a Galerkin projection of the positive definite operator $-c^2(\bx) \Delta$ defined in a bounded and simply connected domain $\Omega \subset \RR^d$. While $\Omega$ may be a real domain, in most applications it should be understood as a truncation of $\RR^d$,  justified using the causality of propagation and the finite duration of the measurements. As long as $\partial \Omega$ is sufficiently far from the source to affect the measurements over the time interval $\cT$, it has no effect on the waveform inversion problem, and we can model it with a convenient homogeneous boundary condition, like Dirichlet. 

The Galerkin projection involves inner products of wave snapshots, evaluated on a uniform time grid with step $\tau$ chosen according to the Nyquist sampling criterion for the probing pulse $f(t)$. We compute these inner products in a data-driven way, using functional calculus on the operator 
$-c^2(\bx) \Delta$ that is self-adjoint in the inner product weighted by $c^{-2}(\bx)$. To avoid the weighting, we use the similarity transformation 
\begin{equation}
p^{(s)}(t,\bx) \mapsto P^{(s)}(t,\bx) = \frac{\bar{c}}{c(\bx)} p^{(s)}(t,\bx),
\label{eq:2.1}
\end{equation}
that acts as an identity near the array, where the wave speed equals the known constant $\bar{c}$. The wave operator 
$\partial_t^2 - c^2(\bx) \Delta$ is transformed by \eqref{eq:2.1} to $\partial_t^2 + A$, where 
\begin{equation}
A = - c(\bx) \Delta \big[ c(\bx) \cdot \big],
\label{eq:2.2}
\end{equation}
is positive definite and self-adjoint, with compact resolvent. 

The array response matrix  \eqref{eq:I3} stores point-wise measurements of the pressure. Because we are interested in inner products, we transform $\bcM(t)$ to  a new data matrix $\bcD(t)$, as stated in \cite[Lemma 1]{borcea2024meets} proved in \cite[Appendix A]{borcea2020reduced} and 
\cite[Section 2.1]{borcea2024meets},
\begin{equation}
\bcD(t) = \bcM^f(t) + \bcM^f(-t), \qquad \bcM^f(t) = - f'(-t) \star \bcM(t).
\label{eq:2.3}
\end{equation}
Here $-f'(-t)$ is the time derivative of the time-reversed probing pulse and $\star$ denotes convolution. The entries of  $\bcD(t)$ have the  inner product expression  \cite[Lemma 1]{borcea2024meets}
\begin{equation}
\mathcal{D}^{(r,s)}(t) = \int_{\Omega} d \bx \, u_0^{(r)}(\bx) u^{(s)}(t,\bx), \qquad 1 \le r,s \le m,
\label{eq:2.4}
\end{equation}
where 
\begin{equation}
u^{(s)}(t,\bx) = \cos \big[t \sqrt{A} \, \big] u^{(s)}_0(\bx),
\label{eq:2.5}
\end{equation}
solves the homogeneous wave equation 
\begin{equation}
\big(\partial_t^2 + A \big)u^{(s)}(t,\bx) = 0, \quad t > 0, ~~ \bx \in \Omega,
\label{eq:2.6}
\end{equation}
with initial conditions
\begin{equation}
u^s(0,\bx) = u_0^{(s)}(\bx) = \big| \hat f \big( \sqrt{A} \, \big) \big| \delta(\bx-\bx_s), \quad \partial_t u^{(s)}(0,\bx) = 0.
\label{eq:2.7}
\end{equation}
In these equations $\hat f$ denotes the Fourier transform of the probing pulse and functional calculus on $A$ is defined using its spectral decomposition: The spectrum of $A$ consists of the 
positive sequence of eigenvalues $(\la_j)_{j \ge 1}$ and the eigenfunctions form an orthonormal basis of $L^2(\Omega)$~\cite[Section 5.3]{kato1966perturbation}. Then, for any continuous function $\psi:\mathbb{R} \mapsto \mathbb{R}$, the  operator $\psi(A)$ is self-adjoint, with the same eigenfunctions as $A$ and eigenvalues $\psi(\la_j)$, for $j \ge 1$.

The Galerkin projection space is spanned by the snapshots of \eqref{eq:2.5} at instants $\{t_j = j \tau, ~ 0 \le j \le n-1\}$, with $n$ chosen according to the depth to which we wish to estimate the medium.  We collect all the snapshots, for all the source excitations, in the $nm-$dimensional row vector field 
\begin{equation}
\bU(\bx) = \big( \bu_0(\bx), \ldots, \bu_{n-1}(\bx) \big), 
\label{eq:2.10}
\end{equation}
where 
\begin{equation}
\bu_j(\bx) = \bu(t_j,\bx) = \big(u^{(1)}(t_j,\bx), \ldots, u^{(m)}(t_j,\bx)\big) \in \RR^{1 \times m}, \qquad 0 \le j \le n-1.
\label{eq:2.8}
\end{equation}
The span of the components of $\bU(\bx)$ is the projection space, written as ``$\mbox{range} \big(\bU(\bx) \big)$". 
The Galerkin approximation of $\bu(t,\bx)$ in this space is 
\begin{equation}
\bu^{\Gal}(t,\bx) = \bU(\bx) \bg(t), 
\label{eq:2.9}
\end{equation}
where the  time-dependent matrix of Galerkin coefficients $\bg(t) \in \RR^{nm \times m}$  satisfies the following system of ordinary differential equations 
\begin{equation}
\int_{\Omega} d \bx \, \bU^T(\bx) \big(\partial_t^2 + A \big) \bu^{\Gal}(t,\bx) = {\bf 0},
 \label{eq:2.11}
\end{equation}
with initial conditions
\begin{equation}
\bu^{\Gal}(0,\bx) = \bu_0(\bx), \quad  \partial_t \bu^{\Gal}(0,\bx) = {\bf 0}.
\label{eq:2.11.i}
\end{equation}
Note that since $c(\bx)$ is unknown, neither $\bU(\bx)$ nor the operator $A$ is known. Nevertheless, it is shown in \cite{Borcea2022rom,borcea2024meets} that the Galerkin system \eqref{eq:2.11} is data-driven and one can obtain an algebraic,  projection ROM
equivalent of the wave equation from it. We state next the definition of the ROM, and then explain its data-driven computation.

Let $\bV(\bx)$ store the causal orthonormal basis of $\mbox{range} \big(\bU(\bx) \big)$. Causality means that $\bV(\bx) = 
\big( \bv_0(\bx), \ldots, \bv_{n-1}(\bx) \big)$ has the same algebraic structure as $\bU(\bx)$,  and 
\begin{equation}
\bv_j(\bx) \in \mbox{span} \big( \bu_0(\bx), \ldots, \bu_j(\bx) \big), \qquad 0 \le j \le n-1.
\label{eq:2.12}
\end{equation}
Such a basis is defined by the Gram-Schmidt orthogonalization, written as 
\begin{equation}
\bU(\bx) = \bV(\bx) \bR,
\label{eq:2.13}
\end{equation}
where $\bR \in \RR^{nm \times nm}$ is block upper triangular. Substituting \eqref{eq:2.13} in 
\eqref{eq:2.11} and multiplying on the left by $\big(\bR^{-1}\big)^T$, we get the 
ROM  version of equation \eqref{eq:2.6}
\begin{equation}
\left( \partial_t^2 + \bA^{\RM} \right) \bu^{\RM}(t) = {\bf 0}.
\label{eq:2.14}
\end{equation}
The true wave $\bu(t)$ is replaced   in this equation by 
\begin{equation}
\bu^{\RM}(t) = \bR \bg(t) = \int_{\Omega} d \bx \, \bV^T(\bx)  \bu^{\Gal}(t,\bx),
\label{eq:2.15}
\end{equation}
and the operator ROM is the $nm \times nm$ matrix 
\begin{equation}
\bA^{\RM} = \int_{\Omega} d \bx \, \bV^T(\bx) A \bV(\bx).
\label{eq:2.16}
\end{equation}

A key point in \cite{Borcea2022rom,borcea2024meets} is that even though the projection space, its basis $\bV(\bx)$ and operator $A$ are unknown, both matrices $\bR$ and $\bA^{\RM}$ are data driven.
Indeed, consider the Gramian of the snapshots, called the ``mass matrix", 
\begin{equation}
\bM = \int_{\Omega} d \bx \, \bU^T(\bx) \bU(\bx) \stackrel{\eqref{eq:2.13}}{=}  \bR^T \underbrace{\int_{\Omega} d \bx \, \bV^T(\bx) \bV(\bx)}_{\mbox{identity}} \, \bR =  \bR^T \bR.
\label{eq:2.17}
\end{equation}
It is an $nm \times nm$ symmetric and positive definite matrix, with Toeplitz+Hankel block structure, and data driven $m \times m$ blocks \cite[Theorem 1]{borcea2024meets}
\begin{equation}
\bM_{i,j} = \int_{\Omega} d \bx \, \bu_i^T(\bx) \bu_j(\bx) = \frac{1}{2} \left[ \bcD(t_{i-j}) + \bcD(t_{i+j}) \right], \quad 0 \le i,j \le n-1.
\label{eq:2.18}
\end{equation}
The block upper triangular matrix 
$\bR$ is the block Cholesky square root of $\bM$, according to equation \eqref{eq:2.17}. 
The $nm \times nm$ ``stiffness matrix" 
\begin{equation}
\bS = \int_{\Omega} d \bx \, \bU^T(\bx) A \bU(\bx) \stackrel{\eqref{eq:2.13},\eqref{eq:2.16}}{=} \bR^T \bA^{\RM} \bR,
\label{eq:2.19}
\end{equation}
is also data driven, with  $m\times m$ blocks determined by \cite[Theorem 2]{borcea2024meets}
\begin{equation}
\bS_{i,j} = \int_{\Omega} d \bx \, \bu_i^T(\bx) A \bu_j(\bx) = - \frac{1}{2} \left[\ddot\bcD(t_{i-j}) + \ddot\bcD(t_{i+j}) \right], \quad 0 \le i,j \le n-1,
\label{eq:2.20}
\end{equation}
where the two dots denote the second derivative. Note from definition \eqref{eq:2.3} that this derivative can be passed in the convolution  to the probing signal. 

The operator ROM is computed in a data driven manner as 
\begin{equation}
\bA^{\RM} = \big(\bR^{-1}\big)^T\bS \bR^{-1},
\label{eq:2.21}
\end{equation}
and the mapping $c(\bx) \stackrel{\cF^{\RM}}{\longrightarrow} \bA^{\RM}$ is deduced from \eqref{eq:2.16}, 
\begin{equation}
\cF^{\RM}[c] = \bA^{\RM} = -\int_{\Omega} d \bx \, \bV^T(\bx) c(\bx) \Delta \big[ c(\bx) \bV(\bx)\big],
\label{eq:2.22}
\end{equation}
where the Laplacian is understood to act columnwise. 

Note how $c(\bx)$ appears explicitly in equation \eqref{eq:2.22}, as the 
coefficient in the operator $A$ defined in \eqref{eq:2.2}, but also implicitly, because the orthonormal basis $\bV(\bx)$ depends on $c(\bx)$.
Our ROMNet approach seeks to learn another map $\cF^{\ML}$, from $\bA^{\RM}$ to another matrix $\bA^{\ML}$,  so that the composition mapping 
\begin{equation}
c(\bx) \stackrel{\cF^{\RM}}{\longrightarrow} \bA^{\RM} \stackrel{\cF^{\ML}}{\longrightarrow} \bA^{\ML}
\end{equation}
gives an explicit, quadratic  dependence of $\bA^{\ML}$ on $c(\bx)$.

\section{The ROMNet approach}
\label{sect:2}
In this section, we introduce ROMNet. We begin with the motivation and the definition of the map $\cF^{\ML}$. Then, we give the details of the algorithm and discuss the advantages and the resolution of our approach.

\subsection{Motivation}
\label{sect:2.1}
We deduce from  
equation \eqref{eq:2.18} that the mass matrix $\bM$ and the data matrices\footnote{Note 
that while $\bU(\bx)$ stores snapshots up to time $(n-1)\tau$, the data go up to the double time $2(n-1)\tau$. This accounts for the round trip travel time, from the array, to points inside the medium and then back to the array.} 
$\{\bcD(t_j), ~ 0 \le j \le 2n-2\}$ are related by a linear mapping that is a bijection onto its range, the set of symmetric Toeplitz+Hankel block matrices. Indeed, \eqref{eq:2.18} gives $\bcD(t_i) = \bM_{i,0}$ for $0 \le i \le n-1$, and the remaining matrices are determined recursively from $\bcD(t_{i+j}) = 2 \bM_{i,j} - \bcD(t_{i-j})$, up to $\bcD(t_{2n-2})$. The matrices $\bcD(t)$ are related to the measurements \eqref{eq:I3}  by the transformation \eqref{eq:2.3}, where  $f'(-t)$ can be deconvolved in principle, because it has the same frequency support as $\bcM(t)$. 
Thus, estimating $c(\bx)$ from $\bM$ is basically the same as the data fitting in FWI. 

We also deduce from equations \eqref{eq:2.2}, \eqref{eq:2.3} and \eqref{eq:2.20}  that the stiffness matrix 
\begin{equation}
\bS  =\int_{\Omega} d \bx \, \bU^T(\bx) A \bU(\bx) = -\int_{\Omega} d \bx \, \bU^T(\bx) c(\bx) \Delta \big[ c(\bx) \bU(\bx)\big],
\label{eq:2.19p}
\end{equation}
is also related to the data  by a linear mapping. It  looks similar to the mapping that gives $\bM$, except for the second time derivative. 
This derivative has no beneficial effect on the estimation of $c(\bx)$, so ``inverting" $c(\bx) \stackrel{\cF^{\Gal}}{\longrightarrow} \bS$ is at least as difficult as FWI. 

Equation \eqref{eq:2.21} shows that the mapping from the data to $\bA^{\RM}$ is nonlinear. While the mass matrix $\bM$ depends linearly on the data, 
  its square root $\bR$ does not, and we also take its inverse. 
The results in \cite{Borcea2022rom,borcea2024meets}
illustrate that the estimation of $c(\bx)$ via the optimization 
\begin{equation}
\mbox{arg}\min_{\hat c \in \mathcal{C}}\,  \mathcal{O}^{\RM}(\, \hat c \, ) + \mbox{regularization}, \qquad  \mathcal{O}^{\RM}(\, \hat c\, ) = \|\cF^{\RM}[\, \hat c\, ] - \bA^{\RM} \|_F^2,
\label{eq:3.1}
\end{equation}
is better than FWI and therefore, better than the ``inversion" of $\cF^{\Gal}$. 

So,  why is  the ``inversion" of $\cF^{\RM}$ better than that of $\cF^{\Gal}$?  We see from equations \eqref{eq:2.19p} and \eqref{eq:2.22} that both $\bS$ and $\bA^{\RM}$ depend explicitly on $c(\bx)$ via the operator $A$ defined in  \eqref{eq:2.2}. However, they differ in their implicit dependence on $c(\bx)$:  via $\bU(\bx)$ in  $\bS$ and via  $\bV(\bx)$ in $\bA^{\RM}$. This is important, because the dependence of $\bU(\bx)$ on $c(\bx)$ appears to be much stronger than that of $\bV(\bx)$. This is a conjecture based on limited theory and extensive numerical simulations.  Here is what is known: 

All the information 
about $c(\bx)$ carried by the measurements of $\bU(\bx)$ is contained in the computable matrix $\bR$.  Once we fix the block Cholesky decomposition algorithm\footnote{The block Cholesky decomposition is unique up to the definition of the diagonal blocks. See \cite[Algorithm 4.1]{untangling} for our implementation.}  to get the square root of $\bM$, we have a bijective mapping between the data and $\bR$.
As seen from equation \eqref{eq:2.17}, the basis $\bV(\bx)$ plays no role in the Gramian and therefore, in the data match.  
%It is stated in \cite[Theorem 1]{borcea2024meets} that the $nm \times m$ block columns of $\bR$ are, in fact, the snapshots in the algebraic, propagator ROM space\footnote{The propagator ROM 
%is not the same as the operator ROM, but the two are related. We refer to \cite[Section 2.4]{borcea2024meets} for the definition of the propagator ROM and its relation to $\bA^{\RM}$.}. The purpose of $\bV(\bx)$ is to transform $\bR$  to the snapshots in the physical space, 
%gathered in $\bU(\bx)$ (recall equation \eqref{eq:2.13}). 

\begin{figure}
\centering
\includegraphics[width=\linewidth]{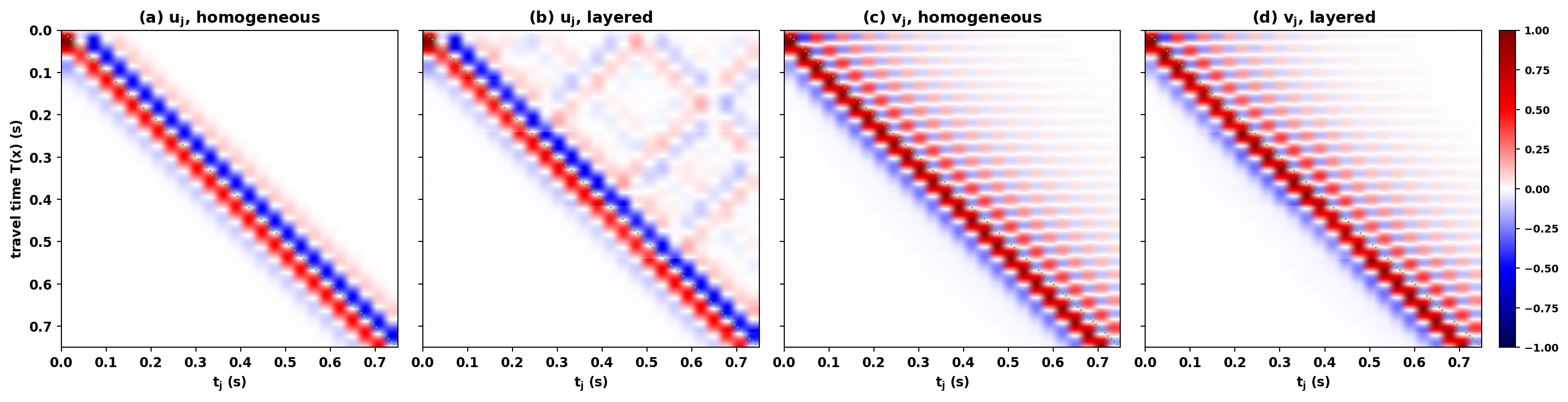}
\vspace{-0.3in}
 \caption{One-dimensional simulations in a homogeneous and layered medium.{The domain is $[0,2.5]$km, with homogeneous Dirichlet boundary conditions, and the source is at $x = 0.06$km. The homogeneous medium has a wave speed $3$km/s. The layered medium has the same wave speed, except in the interval $[0.35, 2.1]$km, which is divided into four equal layers, with wave speeds $3.5, 2.5, 3.5, 2.5$km/s, from top to bottom.} The snapshots in the homogeneous medium (a) and the heterogeneous medium (b). Orthogonalized snapshots in the homogeneous medium (c) and the heterogeneous medium (d).
 The abscissa is the discrete time $t_j = j \tau$, for $j = 0, \ldots, n-1$. The ordinate is the travel time $T$. For each $t_j$, we use colors to plot the amplitudes of $u_j(x(T))$ and $v_j(x(T))$, respectively. 
 }
 \label{fig:1}
 \end{figure}

To gain an intuition about the behavior of $\bV(\bx)$, it is useful to consider the  one dimensional problem, where $m =1$. We give an illustration  in Fig. \ref{fig:1}. The plots Fig. \ref{fig:1}(a)-(b) display the snapshots in a homogeneous medium and in a heterogeneous, layered medium. The abscissa is the discrete time $t_j = j \tau$, for $j = 0, \ldots, n-1$, 
and the ordinate is the travel time $T(x) = \int_0^x c^{-1}(x') dx'$. 
Note the progressing (causal) advancement of the wave front in both media and the reflections in the heterogeneous medium. While the snapshots are very different in the two media, the orthogonal bases plotted in Fig. \ref{fig:1}(c)-(d) are basically the same. To understand why, it helps to consider the sampling of $\bU(\bx)$ on a uniform travel time grid, 
coordinated with the time instants of the snapshots. Such sampling gives an $n \times n$ upper triangular matrix  $\tilde \bU$, whose $j^{\rm th}$ column is the discretization of $u_j(x(T))$. The discrete version of the Gram-Schmidt equation  \eqref{eq:2.13} is the  QR factorization method \cite{golub2013matrix}, that computes an orthogonal matrix $\tilde \bV$,  to  transform the left hand side $\tilde \bU$, to an upper triangular matrix.  But $\tilde \bU$ is already upper triangular, which means that $\tilde \bV$ must be the identity matrix, independent of the reflectivity of the heterogeneous medium. The analysis of $\bV(\bx)$ in the continuum is more involved, but it has been carried out in \cite[Appendix A]{borcea2022reduced},  using the mathematics of waves in layered media \cite{fouque2007wave}.  Explicitly, it is shown there that if the waves are sampled well enough in time i.e., $\tau$ is sufficiently small, 
then  the components of $\bV(x(T))$ are peaked at the wavefront and are  basically the same as those of the 
orthonormal basis computed in a non-scattering medium.  An error estimate of the approximation of $\bV(x(T))$ by the basis in a non-scattering medium  can also be found in  the recent study \cite{druskin2025optimality}.  

\begin{figure}[!htbp]
\centering
\includegraphics[width=\linewidth]{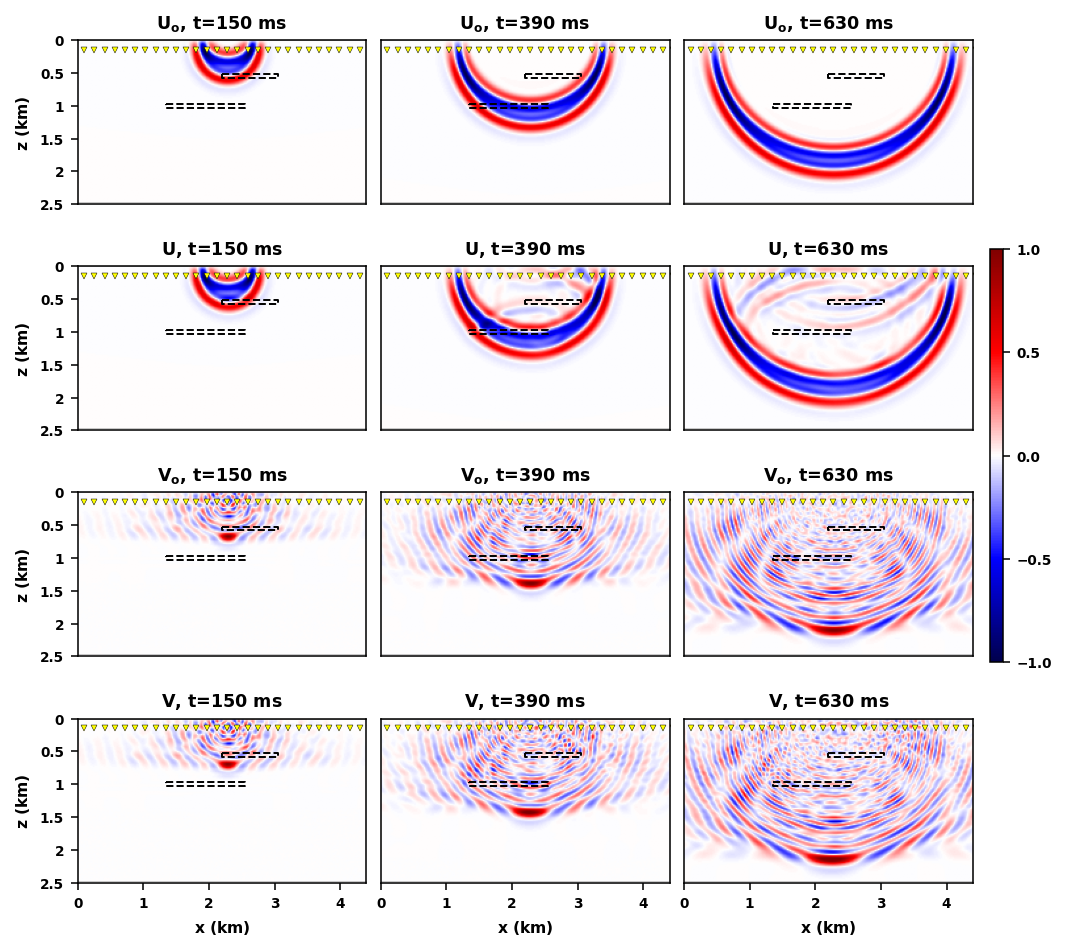}
\vspace{-0.3in}
\caption{Two dimensional simulations. The reference medium is homogeneous, with wave speed $c_o(\bx) \equiv \bar{c} = 3$km/s. The true speed $c(\bx)$ differs from it only inside two thin horizontal bars outlined by the dashed rectangles, where the wave speed is  $4$km/s.  The top two rows show the snapshots  in the reference and heterogeneous media, respectively, at three time instants. The bottom two rows show the corresponding orthonormal basis vectors. The array is linear and lies at the top of the plotted domain, where we have a Dirichlet boundary condition. The subscript $o$ in the titles refers to $c_o$, and the yellow triangles are the sensors in the array. The contour of the reflectors is outlined in all the pictures, including those corresponding to the reference medium,  to aid in the interpretation.}
\label{fig:v-localization}
\end{figure}

The mathematical analysis of $\bV(\bx)$ in higher dimensions remains out of reach, but its dependence on $c(\bx)$ has been investigated numerically \cite{druskin2016direct,borcea2020reduced,Borcea2022rom,borcea2023waveform,borcea2024meets}. We give an illustration of a typical result  in Fig. \ref{fig:v-localization}, for a simple medium with two thin reflectors, where it is easier to visualize the waves. 
The plots in  the first two rows display three snapshots of $u^{(s)}(t,\bx)$ for the center source element in the array. The first row is for the reference medium 
with constant velocity $c_o = \bar{c}$ and shows the circular wavefront moving away from the source. The second row is for the reflective medium, where the snapshots are very different: The wavefront is slightly perturbed, but more importantly, there is a significant back-scattered component. 
The last two rows display the corresponding components from $\bV(\bx)$. They look similar in the two media and display a large peak at the intersection of the 
wavefront with the downgoing ray, starting from the source. They also have weaker oscillations behind the wavefront. These oscillations depend on the 
particular block Cholesky algorithm used to compute $\bR$ i.e., the choice of the orthogonalization of the snapshots corresponding to the same time instant but different sources. We have not explored which choice of the block Cholesky algorithm is better, but in any case, the oscillations have a negligible effect in our study, 
because they arise at scales that are much smaller than the wavelength. We show in section \ref{sect:resolA} that it is impossible to determine features of the medium below the resolution limit of the order of the wavelength. Assuming a parametrization of the search speed $\hat{c}(\bx)$ that respects this resolution limit, the effect of the fast oscillations of the orthonormal bases averages out in the matrix $\bA^{\ML}$ defined below.

The strong similarity of the orthonormal bases in Fig.~\ref{fig:v-localization} is because the two media have similar kinematics. Otherwise, the comparison is less favorable, because the wavefront will be at different locations in the true and reference medium, due to the different kinematics. Still, all numerical evidence suggests that $\bV(\bx)$ changes slowly with the reflectivity, unlike $\bU(\bx)$. 
Furthermore, the components of $\bV(\bx)$ remain peaked at the wavefront, near the forward (down) going rays starting from the sources of excitation. Thus, they give a better localized ``view" of the medium than the snapshots, which are large on the whole transmitted and reflected wavefronts. 

\subsection{Definition of the $\cF^{\ML}$ map}
The observations above suggest that the mapping from the wave speed to $\bV(\bx)$ is smooth and slowly changing. This 
motivates our definition of the map to be learned
\begin{equation}
\bA^{\RM} \stackrel{\cF^{\ML}}{\longrightarrow}  \bA^{\ML} = - \int_{\Omega} d \bx \, \bV(\bx;c_o)^T c(\bx) \Delta \big[ c(\bx) \bV(\bx;c_o)\big],
\label{eq:3.2}
\end{equation}
where $\bV(\bx;c_o)$ is the orthonormal basis computed in a reference medium,  with known wave speed $c_o(\bx)$ that may be constant or variable. Note our notation convention: For the vector fields and matrices computed at the reference or search speed, we indicate that speed in the arguments. For the vector fields and matrices corresponding to the true wave speed $c(\bx)$, we drop the argument, like in   $\bV(\bx;c) \leadsto \bV(\bx)$.

Comparing \eqref{eq:3.2} with \eqref{eq:2.22}, we see that $\bA^{\ML}$ differs from $\bA^{\RM}$ by the fact that it has a simpler, quadratic  dependence on the unknown $c(\bx)$.
Thus, 
we  can estimate $c(\bx)$ from $\bA^{\ML}$ via the optimization, 
\begin{align}
&\mbox{arg}\min_{\hat c \in \mathcal{C}}\,  \mathcal{O}^{\ML}(\, \hat c \, ) + \mbox{regularization}, \nonumber \\
 &\mathcal{O}^{\ML}(\, \hat c\, ) =  \|\cF^{\ML}[\, \bA^{\RM} \, ] +
\int_{\Omega} d \bx \, \bV(\bx;c_o)^T \hat c(\bx) \Delta \big[ \hat c(\bx) \bV(\bx;c_o)\big]\|_F^2.
\label{eq:3.3}
\end{align}

\subsection{The ROMNet algorithm}
\label{sect:MLAlg}
We begin the description of our inversion algorithm with a dimension reduction of the learning problem.  It is proved in \cite[Appendix E]{Borcea2022rom} that the entries of $\bA^{\RM}$ 
decay away from the main diagonal and the proof extends verbatim to $\bA^{\ML}$. The numerical results in  \cite{Borcea2022rom}
show that the estimation of $c(\bx)$ via the optimization \eqref{eq:3.1} does not perform significantly better than the optimization 
that fits  only the first few diagonals of $\bA^{\RM}$. We use this observation to reduce the computational cost, as follows:
Let $W_{\gamma}$ be the  operator that takes an $nm \times nm$ symmetric matrix, 
extracts its entries in the first $\gamma m$ diagonals (the main diagonal and the $\gamma m -1$ super-diagonals) and collects them in a vector of length 
\begin{equation}
N_{\gamma} = 
\sum_{q=0}^{\gamma m -1} (nm - q) = \gamma n m^2 - {\gamma m (\gamma m - 1)}/{2}.
\label{eq:Ngamma}
\end{equation}
We use a neural network with vector $\btheta$ of learned parameters, to approximate the mapping
\begin{equation}
\bA^{\RM}(\hat c) \stackrel{\cG_{\btheta}}{\longrightarrow} W_\gamma\big[\bA^{\ML}(\hat c)\big].
\label{eq:NetLearn}
\end{equation}
The windowing parameter $\gamma$ is user-defined. We illustrate in section \ref{sect:3} the effect of $\gamma$ on the accuracy of the estimated 
wave speed. Most of our simulations are for $\gamma = 5$.

 The neural network is trained on the velocity samples $\{c_i(\bx): i=1,\dots,N_{\rm train} \}$ drawn from the training distribution. These samples are used in the wave equation to simulate the corresponding array data $\bcM(t;c_i)$, from which we can define the data pairs 
$
\bigl\{\, \bigl(\bA^{\RM}(c_i),\,
\bA^{\ML}(c_i)\bigr) : i=1,\dots,N_{\rm train} \,\bigr\}.
$
The operator ROM matrices $\bA^{\RM}(c_i)$ are computed as explained in section \ref{sect:1} and the matrices $\bA^{\ML}(c_i)$ are computed via the formula \eqref{eq:3.2}. {Two different sets of media are drawn from the same distribution: a ``validation set", used to select the parameters $\btheta$ during the training, and a ``held-out test set", which is never used in the training and gives the results reported in section \ref{sect:3}.} The details of the architecture of the neural network and the learning process are in appendix \ref{ap:A.L}.

\begin{figure}[tb]
\centering
\includegraphics[width=1.0\linewidth]{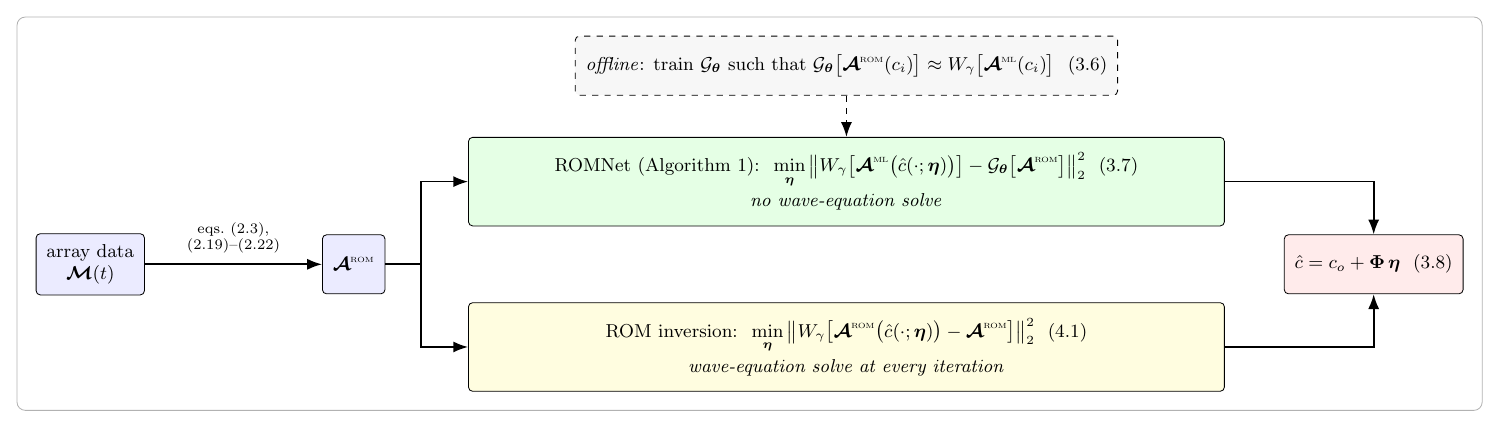}
\caption{ROMNet inversion (top branch) and the  ROM
inversion (bottom branch). Both start from
the data-driven $\bA^{\RM}$. The network $\cG_{\btheta}$ is trained
offline (dashed).
}
\label{fig:pipeline}
\end{figure}
The ROMNet procedure is summarized in Algorithm~\ref{alg:two-stage} and
illustrated in Fig.~\ref{fig:pipeline}. The first step in the algorithm is to compute from the input data $\bcM(t)$ the operator ROM matrix
$\bA^{\RM}$, as explained in section \ref{sect:1}. 
Then, we evaluate the network once to form the Stage~1 prediction {$\cG_{\btheta}\big[\bA^{\RM}\big]$} of {$W_\gamma \big[\bA^{\ML}\big]$} and 
solve via a Gauss-Newton iteration the minimization problem 
\begin{align}
\mbox{arg}\hspace{-0.05in}\min_{\bet \in \RR^{N_c}}  {\mathcal{O}}_{\gamma}(\bet)+ \mbox{regularization}, ~~ {\mathcal{O}}_{\gamma}(\bet) \,=\,
\bigl\|W_{\gamma}\!\bigl[\bA^{\ML}\!\bigl(\hat c(\cdot; \bet) \bigr)\bigr]
- \cG_{\btheta}\bigl[\bA^{\RM}\bigr]\bigr\|_2^{\,2},
\label{eq:Odkhat}
\end{align}
where $\bA^{\ML}\!\bigl(\hat c(\cdot; \bet)\bigr)$ is computed as in equation \eqref{eq:3.2}, with the true $c(\bx)$ replaced by the search speed 
\begin{equation}
\hat c(\bx;\bet) = c_o(\bx) + \boldsymbol{\Phi}(\bx) \bet.
\label{eq:searchC}
\end{equation}
This speed varies about the known reference $c_o(\bx)$, as modeled by a linear superposition of  $N_c$ basis functions gathered in the row vector 
$
\boldsymbol{\Phi}(\bx) = \bigl( \phi_1(\bx), \ldots, \phi_{N_c}(\bx) \bigr),
$
called the ``search dictionary".
The minimization of \eqref{eq:Odkhat} is over the column vector $\bet \in \RR^{N_c}$ of coefficients in the parametrization \eqref{eq:searchC}.

Note that {$\cG_{\btheta}\bigl[\bA^{\RM}\bigr]$} is fixed in \eqref{eq:Odkhat}, so each iteration requires only an evaluation of
{$W_\gamma[\bA^{\ML}(\hat c(\,\cdot\,;\bet))\bigr]$} via the definition \eqref{eq:3.2}. That definition is a Galerkin
projection of the operator $-\hat c(\bx;\bet) \Delta \big[ \hat c(\bx;\bet) \cdot \big]$  on  the fixed basis $\bV(\bx;c_o)$. There is no need to solve 
the wave equation at each iteration, as  in the optimization \eqref{eq:3.1}. 
\begin{algorithm}[!htb]
\caption{ROMNet inversion}
\label{alg:two-stage}
\begin{algorithmic}[1]
\Require array data $\bcM(t)$; trained network $\cG_{\btheta}$; reference
$c_o(\bx)$ and reference basis $\bV(\bx;c_o)$; search dictionary $\bPhi(\bx)$; window parameter $\gamma$;
iteration budget $K$.
\State \textbf{Stage 1: Prediction of $\bA^{\ML}$}
\State \label{ln:computeROM} Compute $\bA^{\RM}$ from $\bcM(t)$.
\State Form the prediction $\cG_{\btheta}\big[\bA^{\RM}\big]$ by a single forward pass through $\cG_{\btheta}$.
\State \textbf{Stage 2: Estimate wave speed}
\State $\bet^{(0)} \gets \mathbf{0}$.
\For{$i=1,\dots,K$}
  \State Form the Gauss--Newton step $\delta\bet^{(i)}$ for problem~\eqref{eq:Odkhat}.
  \State {Backtracking line search:
  $\bet^{(i)} \gets \bet^{(i-1)} + \alpha^\star\delta\bet^{(i)}$.}
\EndFor
\State \Return $\hat c(\bx) = c_o(\bx) + \bPhi(\bx)\,\bet^{(K)}$.
\end{algorithmic}
\end{algorithm}

\subsection{Discussion} 

Although the optimization problem~\eqref{eq:3.3} is simpler than \eqref{eq:3.1}, it is still very difficult to prove that the objective function 
$\mathcal{O}^{\ML}(\, \hat c\, )$ is convex. In general, we remain limited to running numerical simulations to assess the benefit of ROMNet. This is the topic of section \ref{sect:3}. An illustrative, simple assessment is in section \ref{sect:topography}, where we compare two-dimensional cross-sections of the objective functions for the FWI, ROM and ROMNet approaches. We also give in section \ref{sect:resolA} a simple resolution analysis of the estimation of $c(\bx)$ from $\bA^{\ML}$. 

\subsubsection{Landscape of the objective functions}
\label{sect:topography}
\begin{figure}[t]
\centering
\includegraphics[width=\linewidth]{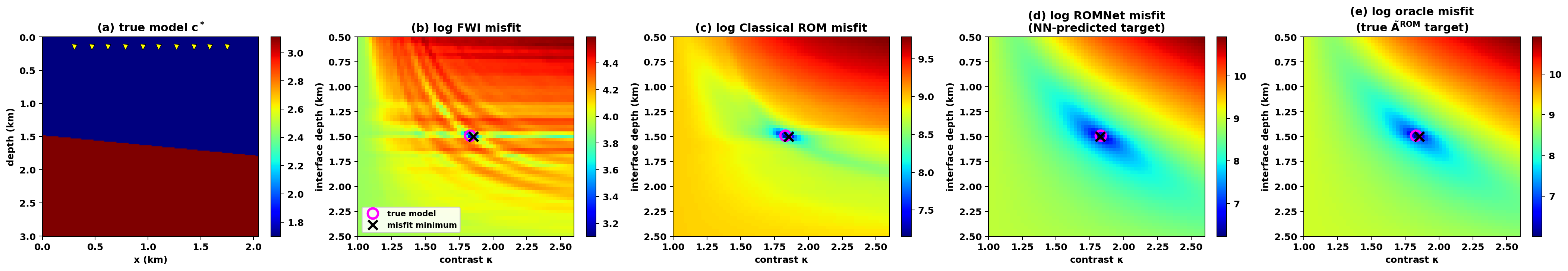}

\vspace{-0.1in} 
\caption{
(a) True model $c(\bx)$ and the linear array located near the top boundary. The abscissa and ordinate are in km and the colorbar is in km/s. (b)--(c) Decimal logarithms of the FWI and ROM objective functions $\mathcal{O}^{\FWI}$, $\mathcal{O}^{\RM}$.
(d)--(e) plot the decimal logarithms of the objective functions $\mathcal{O}^{\ML}$, for $\gamma = 5$: In (d) we use the learned approximation of $\cF^{\ML}$ while in 
(e) we use the exact $\cF^{\ML}$.
The axes in (b)--(e) are: the contrast of the wave speed in the abscissa and the interface position in the ordinate. The global minimum is indicated with  $\textcolor{magenta}{\bigcirc}$ and  $\times$ marks the minimum of each objective on the search grid. The true parameters are not included in the search.
}
\label{fig:landscape}
\end{figure}

To illustrate the landscape of the objective functions, we consider in Fig.  \ref{fig:landscape}(a) a piecewise constant medium. The top region, which hosts the array, has the wave speed $\bar{c} = 1.7$km/s and the bottom region has the wave 
speed $c_{\rm b} = 1.83 \bar{c}$.  The interface position, defined by the depth of the layer at the leftmost point, is $1.487$km. We search 
over the velocity contrast $c_{\rm b}/\bar{c}$ in the interval $[1,2.6]$ and the interface position in the interval $[0.5,2.5]$km. 

The FWI objective function defined in equation \eqref{eq:I4} is displayed in Fig. \ref{fig:landscape}(b). It has multiple (123 on our search grid)  local minima, due to cycle skipping. For the ROM and ROMNet objective functions we use only the $\gamma m$ diagonals of $\bA^{\RM}$ and 
$\bA^{\ML}$, respectively. That is, we display in Fig. \ref{fig:landscape}(c) the ROM objective function
\[
\bigl\| W_\gamma \big[\bA^{\RM}(\hat c) - \bA^{\RM}\big] \bigr\|_2^2 
\]
and  in  Fig. \ref{fig:landscape}(d)--(e)) we show
\[
\bigl\| W_\gamma \big[\bA^{\ML}(\hat c)\big] - \cG_{\btheta}\big[\bA^{\RM}\big]\bigr\|_2^{\,2} \quad \mbox{and} 
\quad \bigl\| W_\gamma \big[\bA^{\ML}(\hat c)-\bA^{\ML}\big] \bigr\|_2^{\,2}.
\]
The last objective function cannot be computed in practice, because $\bA^{\ML}$ depends on the unknown $c(\bx)$. 
We use it as an ``oracle" that gives the best scenario for our method, because by using the true matrix $\bA^{\ML}$, we 
disregard the errors introduced at Stage 1 of Algorithm \ref{alg:two-stage}.

We observe from Fig. \ref{fig:landscape}(c) that the ROM objective function has a better behavior than that of FWI: it has the global minimum at the right 
location and a shallow secondary valley along the diagonal direction. The ROMNet objective function, shown in Fig. \ref{fig:landscape}(d), and the oracle shown in Fig. \ref{fig:landscape}(e), have a  single minimum, the global one. Moreover, the  basin of attraction of the minimum is better than
that of the ROM objective function.

\subsubsection{Resolution analysis} 
\label{sect:resolA}
Even if we manage to estimate the matrix $\bA^{\ML}$ accurately, we need to know how precisely we can determine the wave speed from it,  so we can  use a proper basis in the parametrization \eqref{eq:searchC} of the search speed. 
A typical resolution analysis answers this question in simple media, consisting  of  an isolated inclusion in a homogeneous background,  with wave speed $c_o(\bx) = \bar{c}$.  We give the resolution analysis for the following simple model of an inclusion centered at $\bz$, with radius of support of order $\rho$,
\begin{equation}
c^2(\bx) = \bar{c}^2\Big[ 1 + \big({\sigma}/{(\sqrt{2 \pi} \rho)} \big)^3 \exp\Big(-\frac{|\bx-\bz|^2}{2 \rho^2}\Big) \Big], \quad \bx \in \RR^3.
\label{eq:R1}
\end{equation}
The distance from $\bz$ to the boundary $\partial \Omega$ and the array is assumed
much greater than $\rho$ and the length scale $\sigma$ is used to make the square bracket dimensionless and to control the 
strength of the inclusion. {Note that $\ \int_{\RR^3} \big[ c^2(\bx)/\bar c^2 - 1 \big] d \bx = \sigma^3,$ and that  the Fourier transform of the perturbation is $\sigma^3 \exp(-\rho^2 |\mathbf k|^2/2)$, up to the phase $\exp(-i \mathbf k \cdot \bz)$. This is the source of the two Gaussian factors in Proposition \ref{prop.1}.}

The choice of the Gaussian model \eqref{eq:R1} is convenient for calculations, because if we also consider a modulated Gaussian probing pulse 
\begin{equation}
f(t) = {B}/{\sqrt{2\pi}} \cos(\om_o t) \exp\big[-{(B t)^2}/{2}\big], 
\label{eq:6}
\end{equation}
with central frequency $\om_o$ and bandwidth $B$, then we get a simpler expression of $\bA^{\ML}$. We assume, as is typical in applications, 
that $B \ll \om_o$.

The results are quantitatively similar in any dimension. We carried out the resolution analysis in $\RR^3$,  because it is easier to compute explicitly the matrix $\bA^{\ML}$ using the 
three-dimensional Green's function of the wave equation.  Recall that $\Omega$ is a mathematical truncation of $\RR^3$, whose boundary is not felt during the duration of the measurements. Thus,
we can use the free space Green's function in the analysis.

Using the constant wave speed $c_o(\bx) = \bar{c}$ in definition \eqref{eq:3.2}, we get  that 
\begin{equation}
\bR(\bar{c})^T \bA^{\ML} \bR(\bar{c}) = - \int_{\RR^3} d \bx \, \big[\bU(\bx;\bar{c})\big]^T c(\bx) \Delta \big[ c(\bx) \bU(\bx;\bar{c})\big] = \bGa,
\label{eq:R.3}
\end{equation}
where $\bR(\bar{c})$ is the block Cholesky square root of the mass matrix in the homogeneous medium and $\bU(\bx;\bar{c})$ stores the 
snapshots in the same medium. All the information about  $c(\bx)$ is in  the matrix $\bGa$, whose $(s,r)$ entries in 
the $(j,q)$ block are 
\begin{align}
\Gamma^{(s,r)}_{(j,q)} &= -\int_{\RR^3} d \bx \, u_j^{(s)}(\bx;\bar{c}) c(\bx) \Delta \big[ c(\bx) u_q^{(r)}(\bx;\bar{c})\big],
\label{eq:R4}
\end{align}
for $1 \le s,r \le m$ and $0 \le j, q \le n-1$. We extended the integrals to $\RR^3$, because there are no waves in $\RR^3 \setminus \Omega$. 
For our purpose, it suffices to study the perturbation 
\begin{align}
\delta \Gamma^{(s,r)}_{(j,q)} =& \Gamma^{(s,r)}_{(j,q)}+ \bar{c}^2  \int_{\RR^3} d \bx \, u_j^{(s)}(\bx;\bar{c}) \Delta u_q^{(r)}(\bx;\bar{c}).
\label{eq:R5}
\end{align}
It follows from \cite[Appendix A]{borcea2023waveform} and \cite[Section 2]{druskin2016direct} that 
\begin{equation}
u^{(s)}(t,\bx;\bar{c}) =  f(t) \star \left[G(t,\bx,\bx_s) + G(-t,\bx,\bx_s)\right],
\label{eq:8.1}
\end{equation}
where  $G(t,\bx,\bx_s)$ is the Causal  Green's function, satisfying  
\begin{align}
\left[\partial_t^2 - \bar{c}^2 \Delta\right]G(t,\bx,\bx_s) &= \delta'(t) \delta(\bx-\bx_s), \qquad t \in \RR, ~~ \bx \in \RR^3,\label{eq:G1} \\
G(t,\bx,\bx_s) &\equiv 0, \qquad t < 0. \label{eq:G1.ini}
\end{align}
The second term in \eqref{eq:8.1}, with the $-t$ argument, is due to the even extension in time that is used to define our ROM. We are interested 
in evaluating the wave at time $t > t_f$. The second term in \eqref{eq:8.1} does not contribute at such time,  due to the initial condition \eqref{eq:G1.ini}, so using the known expression of the three-dimensional Green's function we get, for $t>0$ outside the support $(-t_f,t_f)$ of $f(t)$, 
\begin{align}
u^{(s)}(t,\bx;\bar{c}) &= f(t) \star_t G(t,\bx,\bx_s) = 
\frac{f' \big(t - |\bx-\bx_s|/\bar{c}\big)}{4 \pi{\bar c^2} |\bx-\bx_s|}.
\label{eq:8}
\end{align}

Substituting the expression \eqref{eq:8} and the wave speed \eqref{eq:R1} in equation \eqref{eq:R5}, we obtain after a long but straightforward calculation, that involves 
multiple Gaussian integrals and manipulation of Bessel functions, the result stated in the next proposition. This result is an approximation 
of the matrix $\delta \bGa$, obtained in an asymptotic scaling regime that is typical in applications.  This regime is described by the ordering  $L  \gg a \gg \la_o$
of three fundamental length scales: 
The central wavelength $\la_o = 2 \pi \bar{c}/\om_o$, the aperture $a$ of the array (assumed planar and square) and the range of propagation $L \sim |\bx_s-\bz|$. 
The inclusion is strong i.e.,  
$
\sigma/\rho \gg 1,
$
to have  a significant effect, and its radius $\rho$ can be smaller than $\la_o$, of the same order as $\la_o$ or larger than $\la_o$, as long as  it satisfies 
$
\rho \ll \min\{\sqrt{\la_o L},\bar{c}/B\},
$
so that we have simpler formulas. 

\vspace{0.05in}
\begin{prop}
\label{prop.1}
In the asymptotic scaling regime described above, and for the wave speed model \eqref{eq:R1}, we have 
\begin{align}
\delta \Gamma^{(s,r)}_{(j,q)}&\approx \frac{\mathcal{K} \sigma^3}{|\bz-\bx_s| |\bz-\bx_r|} \exp \left\{\hspace{-0.03in}- \frac{B^2}{2} \Big[ \Big(t_j - \frac{|\bz-\bx_s|}{\bar{c}}\Big)^2 + \Big(t_q - \frac{|\bz-\bx_r|}{\bar{c}}\Big)^2\Big] \right\}\nonumber \\
&\hspace{-0.4in}\times \Big\{\Big[ \frac{3}{\rho^2} + \frac{16 \pi^2 \cos^2\big( \varphi_{r,s}/2\big)}{\la_o^2} \Big]\cos\Big\{\frac{2\pi}{\la_o} \Big[\bar{c}(t_j -t_q) -|\bz-\bx_s|+|\bz-\bx_r|\Big]\Big\}  \nonumber \\
&\hspace{1.9in}\times \exp\Big[-8\pi^2  \Big(\frac{\rho\sin(\varphi_{r,s}/2)}{\la_o}\Big)^2\Big] \nonumber \\
&\hspace{-0.4in}- \Big[ \frac{3}{\rho^2} +\frac{16 \pi^2\big[\cos^2\varphi_{r,s}-\cos^2\big(\varphi_{r,s}/2\big)\big] }{\la_o^2} \Big] 
\cos\Big\{\frac{2 \pi}{\la_o} \Big[\bar{c}(t_j +t_q) -|\bz-\bx_s|- |\bz-\bx_r|\Big]\Big\} \nonumber \\
&\hspace{1.9in}\times \exp\Big[-8 \pi^2\Big(\frac{\rho \cos(\varphi_{r,s}/2)}{\la_o}\Big)^2\Big]
\Big\}, \label{eq:result}
\end{align}
where $\varphi_{r,s}$ denotes the angle
between the vectors $\bz-\bx_s$ and $\bz-\bx_r$ and 
$
\mathcal{K} = \Big(\frac{\bar{c} \om_o B}{16 \pi^{3/2}} \Big)^2.
$
\end{prop}

\vspace{0.1in}
We conclude from the proposition that due to the first exponential, the entries $\delta \Gamma^{(s,r)}_{(j,q)}$ are negligible unless 
\begin{equation}
|\bz-\bx_s| = \bar{c} t_j + O(\bar{c}/B) \quad \mbox{and}  \quad |\bz-\bx_r| = \bar{c} t_q + O(\bar{c}/B). 
\label{eq:Res1}
\end{equation}
The scale $\bar{c}/B$ is the typical range resolution encountered in imaging with reverse time migration or matched filtering  \cite{biondi20063d,jakowatz2012spotlight,cheney2009fundamentals}. 
For times $t_j$ and $t_q$ satisfying \eqref{eq:Res1}, we have three distinct regimes: 

$\bullet$ If the inclusion is very small, with radius
$
\rho \ll \la_o,
$
then  the result simplifies to 
\begin{align*}
\delta \Gamma^{(s,r)}_{(j,q)}\approx \frac{6 \mathcal{K} \sigma^3/\rho^2}{|\bz-\bx_s| |\bz-\bx_r|} &\exp \Big\{ - \frac{B^2}{2} \Big[ \big(t_j - |\bz-\bx_s|/\bar{c}\big)^2 + \big(t_q - |\bz-\bx_r|/\bar{c}\big)^2\Big]\Big\} \nonumber \\
&\times \sin \Big[\frac{2 \pi}{\la_o} \Big(\bar{c} t_j -|\bz-\bx_s|\Big)\Big]\sin \Big[\frac{2 \pi}{\la_o} \Big(\bar{c} t_q -|\bz-\bx_r|\Big)\Big].
\end{align*}
The two sine functions determine $\bz$ with resolution $O(\la_o)$, if we consider two different source-receiver pairs, but $\sigma$ and $\rho$ cannot be determined
separately: only the combination $\sigma^3/\rho^2$ enters. Therefore, the strength and radius of very small inclusions cannot be identified unambiguously.

$\bullet$ If the inclusion is large,  with radius  $\rho \gg \la_o$,  we note that the right hand side in \eqref{eq:result} is significant only for nearby source and receiver pairs,
 that give 
\[
\sin(\varphi_{r,s}/2) = O\Big( \frac{|\bx_s-\bx_r|}{|\bz-\bx_s|} \Big) = O \Big(\frac{\la_o}{\rho}\Big).
\]
The rate of exponential decay in the angle $\varphi_{r,s}$ can be used, in principle, to determine $\rho$. 
On the diagonal $(s=r)$ we have 
\begin{align*}
\delta \Gamma^{(s,s)}_{(j,q)}\approx &\frac{16 \pi^2\mathcal{K} \sigma^3\cos\big[\om_o (t_j -t_q)\big]}{|\bz-\bx_s|^2\la_o^2} \exp \Big\{ \hspace{-0.05in}- \frac{B^2}{2} \Big[ \Big(t_j - \frac{|\bz-\bx_s|}{\bar{c}}\Big)^2 \hspace{-0.02in}+ \Big(t_q - \frac{|\bz-{\bx_s}|}{\bar{c}}\Big)^2\Big]\Big\}
\end{align*}
and only the $(j,q)$ blocks corresponding to
\[
|t_j-t_q| = O(1/B), \quad t_j = |\bz-\bx_s|/\bar{c} + O(1/B),
\]
contribute. The parameter $\sigma$ can be determined from the peak value of $  \delta \Gamma^{(s,s)}_{(j,j)}$ and a better estimation of $\bz$, 
with resolution $O(\la_o)$, comes from the entries corresponding to nearby source receiver pairs, like $(s,r=s\pm 1)$, by exploiting the oscillations of the first cosine in \eqref{eq:result}.

$\bullet$ If $\rho \sim \la_o$, all the terms in   \eqref{eq:result}
contribute, and $\sigma$, $\rho$ and $\bz$ can be determined. 

\vspace{0.1in}
In conclusion, we cannot determine features of $c(\bx)$ at micro scales $\ll \la_o$, but it is feasible to estimate 
variations on scales $\gtrsim \la_o$. This is consistent with the well-known Abbe-Rayleigh resolution limit and suggests a 
parametrization of $\hat c(\bx)$ using basis functions with support $O(\la_o)$.  
In the numerical simulations, the search dictionary has basis functions of smaller support. This allows a more flexible 
representation of the medium, but since we over-parametrize $\hat c(\bx;\bet)$, the Gauss-Newton iteration is regularized, as explained in appendix \ref{sect:GaussN}.

\section{Numerical results}
\label{sect:3}
In this section we present numerical results that assess the performance of ROMNet. 
The assessment includes a comparison with the operator ROM inversion method introduced in \cite{Borcea2022rom}, 
that solves the minimization problem 
\begin{align}
\mbox{arg}\hspace{-0.05in}\min_{\bet \in \RR^{N_c}}  {\mathcal{O}}^{\RM}_{\gamma}(\bet)+ \mbox{regularization}, ~~ {\mathcal{O}}^{\RM}_{\gamma}(\bet) \,=\,
\bigl\|W_{\gamma}\!\bigl[\bA^{\RM}\!\bigl(\hat c(\cdot; \bet) \bigr)
- \bA^{\RM}\bigr]\bigr\|_2^{\,2}.
\label{eq:gammaROM}
\end{align}
We also compare the results with those given by two representative deep learning approaches to FWI: \emph{Fourier-DeepONet} \cite{zhu2023fourier},  a
DeepONet operator network with Fourier-neural-operator layers in its
decoder; and 
\emph{InversionNet} \cite{Wu2019inversionnet,Deng2022openfwi}, a strided convolutional encoder--decoder that is the standard
end-to-end FWI baseline in the OpenFWI benchmark. These use the same training velocity sets as ROMNet, but are trained and evaluated with the publicly released code of the respective publications.

To assess how far the ROMNet estimates are from the best possible ones,  we also show the result of the oracle (recall section \ref{sect:topography}), where we minimize the objective function that is similar to~\eqref{eq:Odkhat}, but the neural network output $\cG_{\btheta}\big[\bA^{\RM}\big]$ is replaced by {$W_\gamma\bigl[\bA^{\ML}\big]$ computed from equation \eqref{eq:3.2}}, at the true $c(\bx)$. The oracle cannot be used in practice,
because $c(\bx)$ is unknown. Nevertheless, it is useful to see how it performs, because it removes the error at Stage 1 of the Algorithm \ref{alg:two-stage}, and thus quantifies the performance 
of the Gauss-Newton method at Stage 2.

{We refer to appendix \ref{sec:exp-setup} for the details of the  wave equation solver and the computation 
of the array response matrix $\bcM(t)$. The neural network and training are described in appendix \ref{ap:A.L}. The Gauss-Newton iteration for solving the minimization problems 
\eqref{eq:Odkhat} and \eqref{eq:gammaROM} is summarized in appendix \ref{sect:GaussN}. The iteration   is run for {$K = 40$ steps for the ROM inversion problem \eqref{eq:gammaROM} and $K = 10$ steps for the ROMNet optimization \eqref{eq:Odkhat}. This suffices because the error changes by less than $1\%$ between the tenth and the fortieth iteration (Fig.~\ref{fig:rg-speed-rel-convergence-p785}). There is a significant computational cost difference per iteration for the two problems, as explained in section \ref{sec:exp-speed}. This difference stems from the fact that  the evaluation of the 
objective function in \eqref{eq:gammaROM} requires solving the wave equation at the current search speed $\hat{c}(\bx;\bet)$, whereas the objective function in \eqref{eq:Odkhat} can be computed directly from formula \eqref{eq:3.2}.

The search speed is parametrized as in equation \eqref{eq:searchC}. 
We use anisotropic Gaussian functions in our simulations, as described in Table \ref{tab:setup-inversion} in appendix \ref{sec:exp-setup}.}
The results are  similar for other search dictionaries (see section \ref{app:hat} for an illustration).

\subsection{Random Gaussians}
\label{sec:exp-rg}
The first training set is given by a collection of wave speeds $\{c_i(\bx), ~i = 1, \ldots, N_{\rm train}\}$
modeled as a random superposition of isotropic Gaussian functions added to the constant reference 
speed $c_o(\bx) = \bar{c} = 3$km/s. 

We assign a Gaussian to each {point of the uniform finite-difference grid of the wave solver (Table~\ref{tab:setup-rg})}, with amplitude drawn uniformly and independently
from the interval $[-4,4]$km/s, and standard deviation chosen randomly, with equal probability, 
from the set of $5$ values that are uniformly log-spaced
between $\lambda_{\rm min}/9\approx 0.033$\,km and $\lambda_{\rm min}=0.3$\,km. Here 
$\lambda_{\rm min} = 2 \pi \bar{c}/(\om_o + B)$ is the smallest wavelength for the frequency content of our pulse.

Because the medium near the array and the boundary is supposed to have a constant wave speed $\bar{c}$, we taper the amplitudes smoothly to zero within $0.375$\,km of the top and bottom boundaries and $0.456$\,km of the lateral boundaries. The sum of the Gaussians is added to $c_o$ and the result is
clipped through a hard cutoff (which is rarely active) to $[1.6,4.0]$\,km/s. Two examples of the resulting wave speeds are shown in the first column of Fig. ~\ref{fig:indist-panels}.

\subsubsection{In-distribution evaluation on the test set}
\label{sec:exp-indist}

Table~\ref{tab:indist-per-sample} summarizes the {mean $L^2$ relative error} of the estimated wave speed, for $11$ media drawn at random from the held-out test set. We also report in Table ~\ref{tab:indist-metrics-appendix} three additional 
error metrics for the learning-based methods. These are used by the OpenFWI benchmark
\cite{Deng2022openfwi} and give the mean absolute error (MAE), the root-mean-square error
(RMSE), and the structural similarity index (SSIM, computed on the interior
with the default $7\times7$ window). ROMNet achieves the smallest error of the three learning-based methods.
\begin{table}[!htp]
\centering
\small
\begin{tabular}{lc}
\toprule
method & relative $L^2$ error \\
\midrule
ROMNet       & \textbf{0.0297} \\
Fourier-DeepONet        & 0.0332 \\
InversionNet            & 0.0419 \\
\midrule
 ROM           & 0.0293 \\
Oracle & 0.0232 \\
\bottomrule
\end{tabular}
\caption{Mean relative $L^2$ error for all methods, over11 held-out test samples. }
\label{tab:indist-per-sample}
\end{table}
\begin{table}[!htp]
\centering
\small
\begin{tabular}{lcccc}
\toprule
method & rel.\ $L^2$ & MAE & RMSE & SSIM \\
\midrule
ROMNet       & \textbf{0.0297} & \textbf{0.0476} & \textbf{0.0890} & \textbf{0.7443} \\
Fourier-DeepONet        & 0.0332 & 0.0517 & 0.0995 & 0.7438 \\
InversionNet            & 0.0419 & 0.0676 & 0.1258 & 0.5902 \\
%\midrule
\bottomrule
\end{tabular}
\caption{Complementary error metrics over the 11 held-out test samples. }
\label{tab:indist-metrics-appendix}
\end{table}

Fig.~\ref{fig:indist-panels} shows the recovered wave speed  on
two representative held-out test samples. ROMNet performs close to the ROM inversion and to the Oracle, and better than the other learning-based methods.

\begin{figure}[t]
\centering
\includegraphics[width=\linewidth]{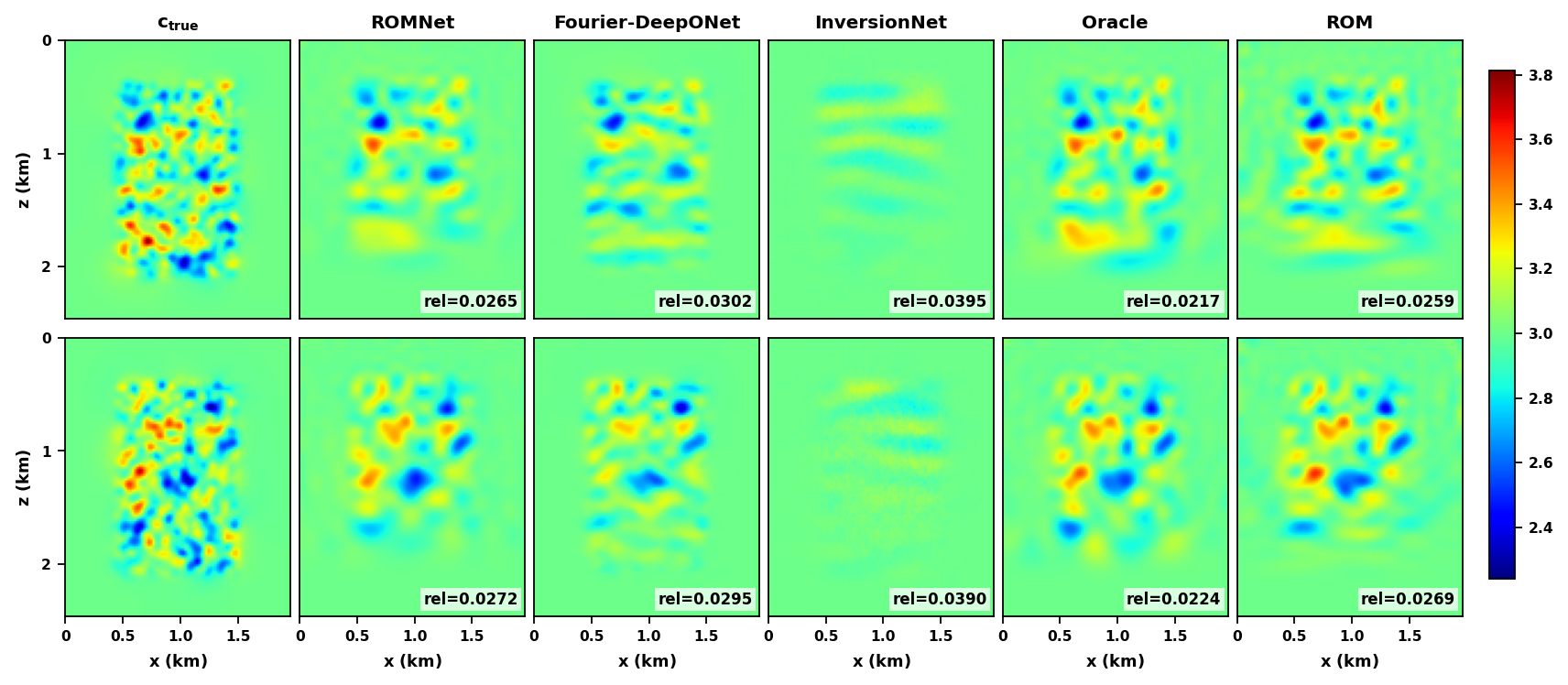}

\vspace{-0.1in}
\caption{Recovered wave speed  on two held-out Random-Gaussian
test samples. From left to right, the columns show ground truth $c_{\rm true}$, ROMNet, Fourier-DeepONet, InversionNet, Oracle, and the  ROM
inversion results. {The $L^2$ relative error of each estimate is printed in its panel.}}
\label{fig:indist-panels}
\end{figure}

\vspace{-0.1in}\subsubsection{Out-of-distribution evaluation}
\label{sec:exp-ood}

We probe the out-of-distribution behavior on piecewise constant media, with three
inclusion shapes: a disk, known in the Geophysics literature as the ``Camembert model" \cite{gauthier1986two}, 
a vertical bar, and a horizontal bar. The 
results are  in Fig. \ref{fig:ood-panels}. Again, ROMNet  attains the lowest error {among the learning-based methods} on each probe. The 
reconstruction is not as accurate as that of the ROM inversion and the Oracle, but it does identify well the inclusions. 

\begin{figure}[t]
\centering
\includegraphics[width=\linewidth]{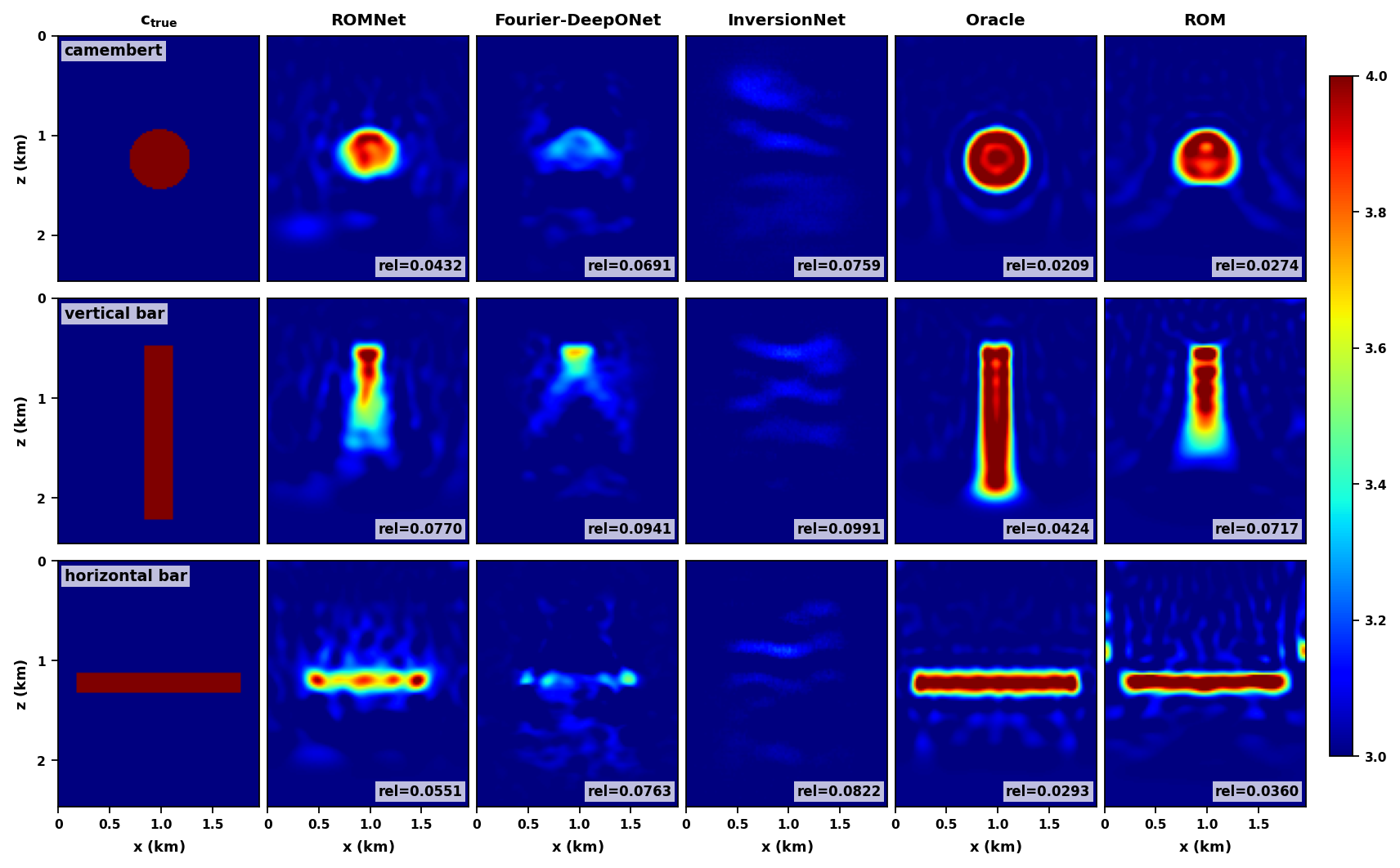}

\vspace{-0.1in}\caption{Recovered velocity on three out-of-distribution probes. From left to right: true 
$c(\bx)$, ROMNet, Fourier-DeepONet, InversionNet,
Oracle, and the  ROM inversion results. {The $L^2$ relative error of each estimate is printed in its panel.}}
\label{fig:ood-panels}
\end{figure}

\subsubsection{Effect of the windowing parameter}
\label{sec:exp-gamma}
As we mentioned already, the numerical results are for the windowing parameter $\gamma = 5$. 
To motivate this choice, we show in Table \ref{tab:dband} the mean $L^2$  relative error for 
$1 \le \gamma \le 5$. The error saturates at $\gamma = 3$ for ROMNet. 
The dependence of the  ROM inversion on $\gamma$ was tested already in \cite{Borcea2022rom},
and there is no significant advantage of taking $\gamma > 5$.
\begin{table}[h]
\centering
\small
\begin{tabular}{lccccc}
\toprule
 & $\gamma=1$ & $\gamma=2$ & $\gamma=3$ & $\gamma=4$ & $\gamma=5$ \\
\midrule
ROMNet & 0.0304 & 0.0299 & 0.0298 & 0.0297 & 0.0297 \\
Oracle & 0.0241 & 0.0233 & 0.0231 & 0.0231 & 0.0232 \\
 ROM & 0.0301 & 0.0293 & 0.0292 & 0.0291 & 0.0289 \\
\bottomrule
\end{tabular}
\caption{Effect of the windowing parameter $\gamma$ on the mean $L^2$  relative error of the reconstruction. {The ROMNet and Oracle rows are averaged over the same $11$ held-out test media as Table~\ref{tab:indist-per-sample}, and the ROM row over $3$ of them, which is why its value at $\gamma = 5$ differs slightly from the one in Table~\ref{tab:indist-per-sample}.}}
\label{tab:dband}
\end{table}

\vspace{-0.1in}
\subsubsection{Effect of the parametrization basis}
\label{app:hat}
A natural alternative to the Gaussian basis functions in $\Phi(\bx)$ is the piecewise-bilinear
``hat'' basis. We display here the results with  such a basis, for the same grid on which we 
centered the Gaussian basis functions. This is a less smooth representation of $\hat c(\bx;\bet)$, 
that one may hope will give a sharper estimation of the medium. However, the reconstructions are 
qualitatively similar, as illustrated in Fig. \ref{fig:hat-rg}.
\begin{figure}[h]
\centering
\includegraphics[width=\linewidth]{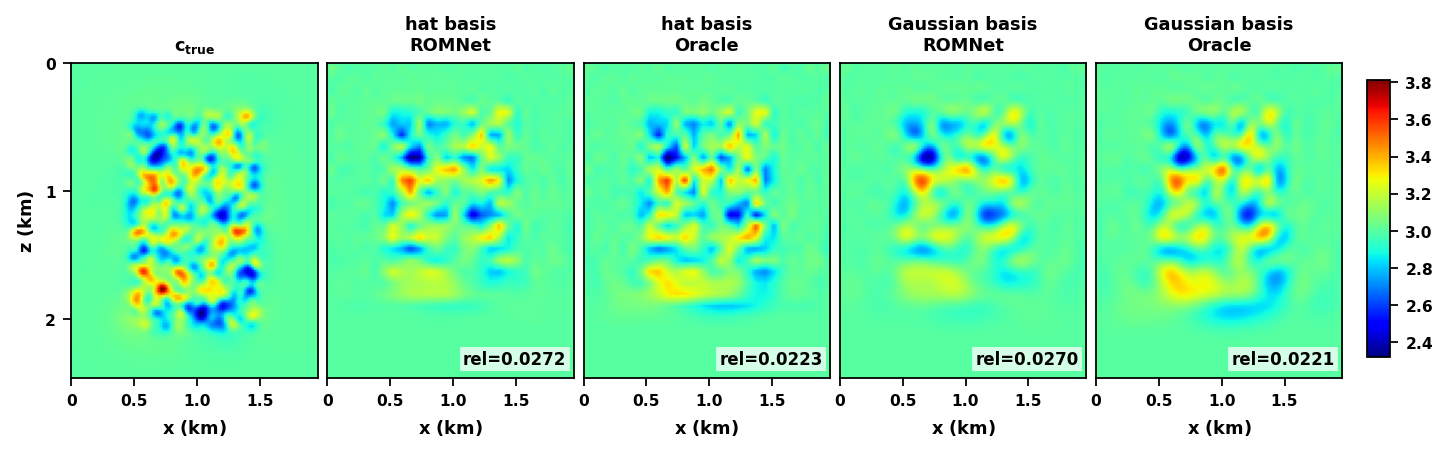}

\vspace{-0.1in}\caption{From
left to right: true  $c(\bx)$, the hat-basis inversion with ROMNet and the Oracle; the Gaussian-dictionary
inversion with ROMNet and the Oracle.}
\label{fig:hat-rg}
\end{figure}

\subsection{GeoFWI}
\label{sec:exp-geofwi}
The second training data set, called GeoFWI, was introduced in \cite{li2026geofwi} to
complement the widely used OpenFWI benchmark \cite{Deng2022openfwi}.  GeoFWI consists of 
wave speed models generated from geological sedimentary and structural
rules, rather than from transformations of natural images. It includes geological
structures like folds, faults, and salt bodies with stacked layers. Such structures were absent in the  random Gaussians data set considered in section \ref{sec:exp-rg}.

We used a ``shallow-water" subset of  GeoFWI, with wave speeds satisfying 
\begin{equation}
c_i(\bx) \,<\, 2.0\,\mathrm{km/s}\quad \mbox{for} ~ \bx = (x,y)  ~\mbox{and} ~y < \la_{\rm min}.
\label{eq:water-shallow}
\end{equation}
Recall from Table \ref{tab:setup-geofwi} that the array lies at depth $\la_{\rm min}$. The condition \eqref{eq:water-shallow}
keeps the wave speed close to that in water, hence the name of the training subset. The reference speed $c_o(\bx)$ varies only in the depth coordinate $y$ of $\bx = (x,y)$, and is defined as the per-depth median of the training velocities in the shallow-water training set.

Table~\ref{tab:geofwi-water-shallow} reports the mean $L^2$
relative error of the estimated speed, evaluated on $45$ samples, $15$ drawn at random from each of the  three categories of wave speed candidates:  nearly layered media, media with faults and media with a salt body. Fig.~\ref{fig:geofwi-water-shallow-folds}--\ref{fig:geofwi-water-shallow-salts}
show the recovered wave speed for three sample speeds drawn from each category. The  ROM inversion starts with an initial wave speed given by 
the smoothed true $c(\bx)$, obtained by convolving it with a Gaussian of standard deviation $133$m in depth and $17$m in the horizontal direction and then averaging in the horizontal direction. This initial guess is different for each medium and uses the unknown $c(\bx)$, whereas ROMNet starts from the reference $c_o(\bx)$, which is computed once from the training set and is the same for all the media. See Fig.~\ref{fig:geofwi-inits} for an illustration. 

Again, the error is smaller for ROMNet than all other 
learning-based methods. {ROMNet and Fourier-DeepONet give comparable reconstructions for the first two categories 
(Fig. \ref{fig:geofwi-water-shallow-folds}--\ref{fig:geofwi-water-shallow-faults}), while InversionNet is clearly worse. For the media with a salt body, ROMNet performs significantly better than the other learning-based methods 
(Fig. \ref{fig:geofwi-water-shallow-salts})}.

\begin{table}[h]
\centering
\small
\setlength{\tabcolsep}{4.5pt}
\begin{tabular}{lcccccc}
\toprule
category & sample size & \shortstack{ROMNet} & \shortstack{Fourier-\\DeepONet} & InversionNet & \shortstack{\\ROM} & Oracle \\
\midrule
fold  & 15 & \textbf{0.0264} & 0.0284 & 0.0528 & 0.0235 & 0.0180 \\
fault & 15 & \textbf{0.0382} & 0.0459 & 0.0799 & 0.0252 & 0.0186 \\
salt  & 15 & \textbf{0.0428} & 0.0523 & 0.0967 & 0.0530 & 0.0241 \\
\midrule
overall & 45 & \textbf{0.0358} & 0.0422 & 0.0765 & 0.0339 & 0.0202 \\
\bottomrule
\end{tabular}
\caption{Mean $L^2$ relative error on $15$ random samples per structural
category of the GeoFWI test set ($45$ in total). Best learning-based method
in each row is in bold. }
\label{tab:geofwi-water-shallow}
\end{table}

\begin{figure}[!htp]
\centering
\includegraphics[width=0.95\linewidth]{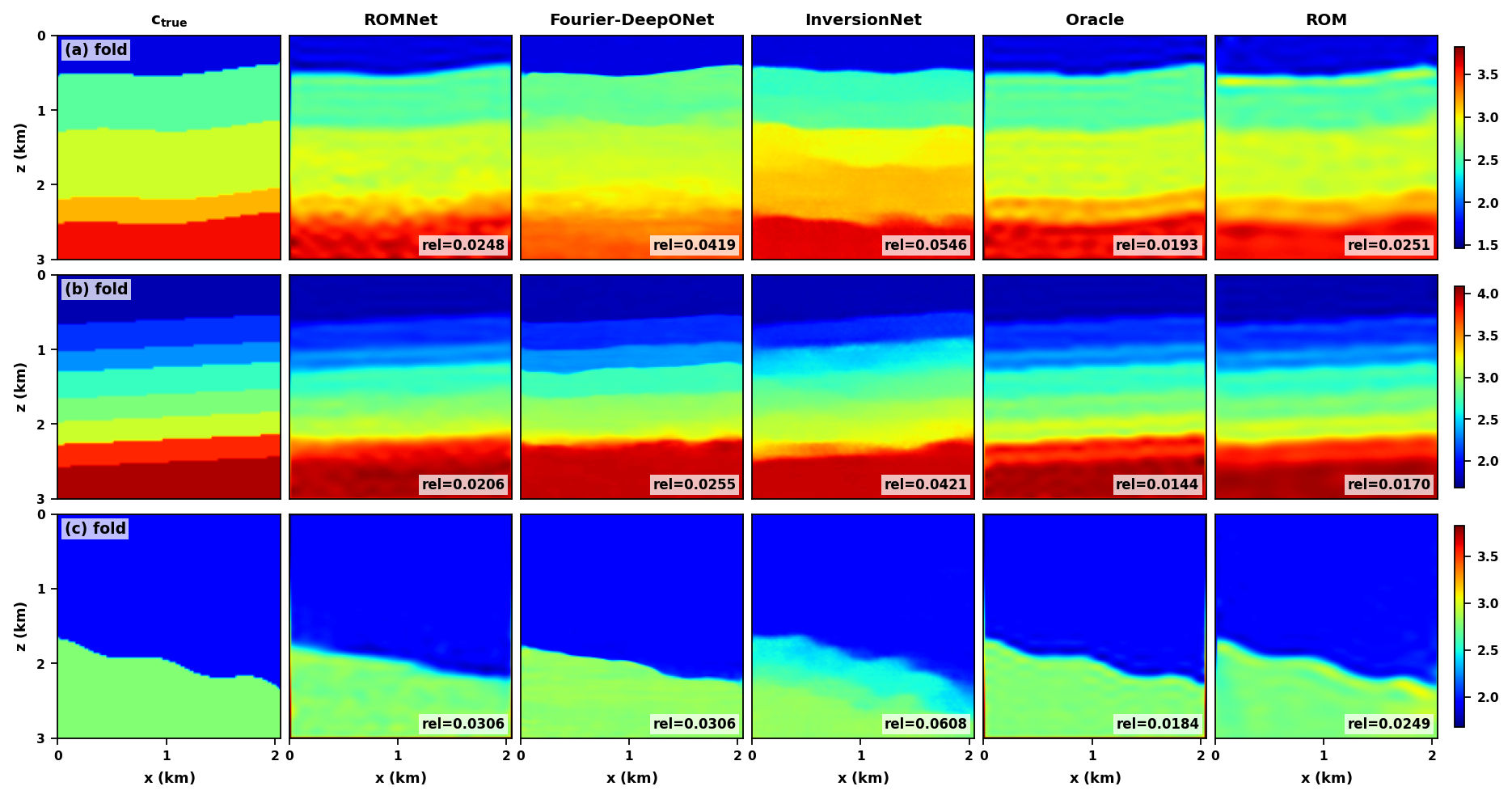}
\vspace{-0.1in}
\caption{Three velocity models with folded layering.
From left to right: True  $c(\bx)$, ROMNet, Fourier-DeepONet, InversionNet,
Oracle and ROM inversion.}
\label{fig:geofwi-water-shallow-folds}
\end{figure}

\begin{figure}[!htp]
\centering
\includegraphics[width=0.95\linewidth]{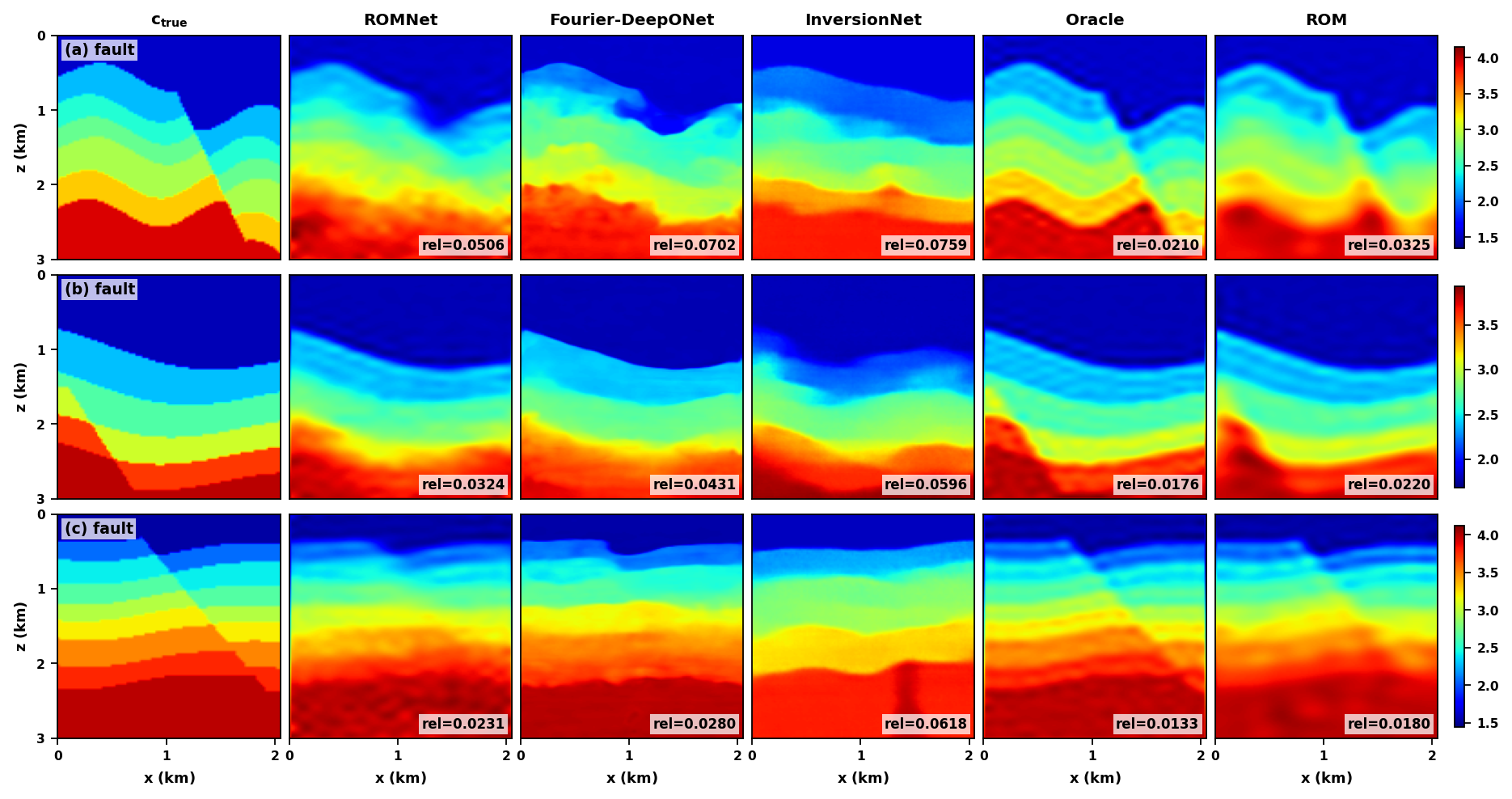}

\vspace{-0.1in}\caption{Three velocity models with faulted layering.
Columns as in Figure~\ref{fig:geofwi-water-shallow-folds}.}
\label{fig:geofwi-water-shallow-faults}
\end{figure}

\begin{figure}[!htp]
\centering
\includegraphics[width=0.95\linewidth]{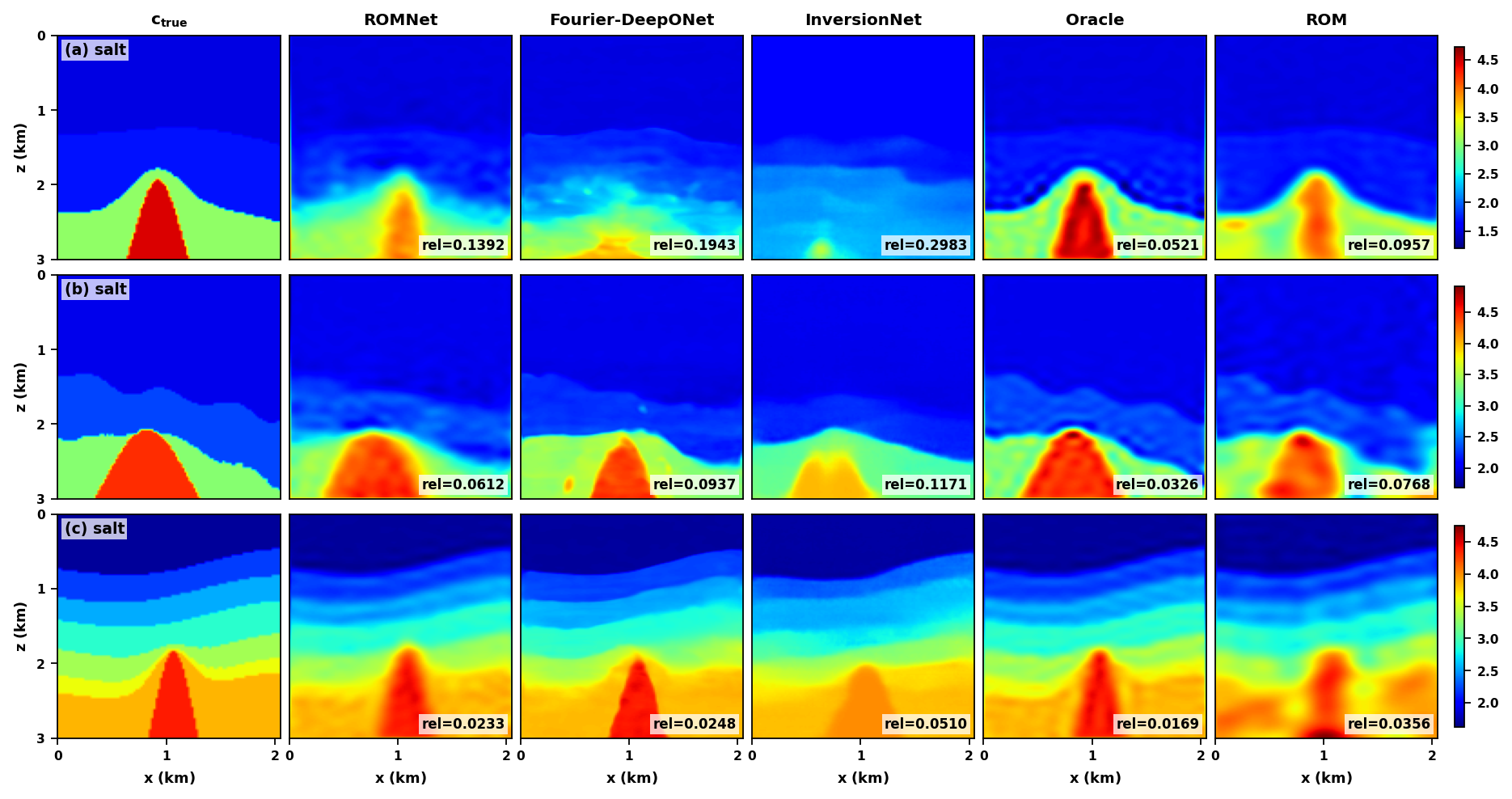}

\vspace{-0.1in}\caption{Three velocity models containing a salt body.
 Columns as in Figure~\ref{fig:geofwi-water-shallow-folds}.}
\label{fig:geofwi-water-shallow-salts}
\end{figure}

\begin{figure}[!htp]
\centering
\includegraphics[width=0.5\linewidth]{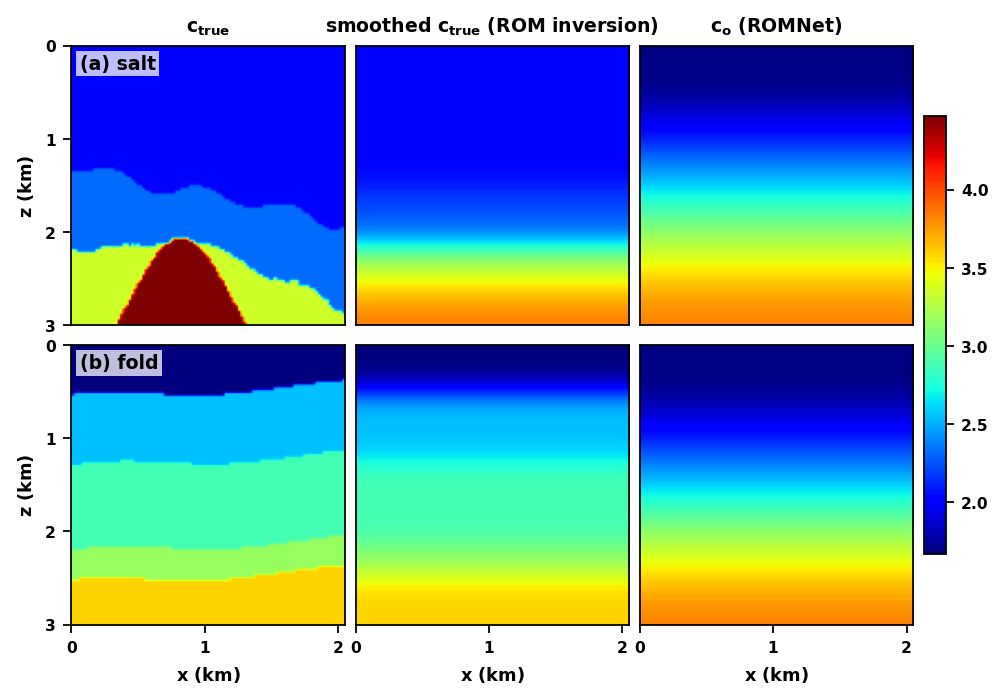}

\vspace{-0.1in}
\caption{GeoFWI: the initial wave speed of the ROM inversion (smoothed true $c(\bx)$, middle column) and the reference $c_o(\bx)$ of ROMNet (right column), for a medium with a salt body (top) and a folded medium (bottom). The true wave speeds are in the left column. Both initial wave speeds depend only on depth; the ROM inversion uses a different one for each medium, whereas $c_o(\bx)$ is the same for all media.}
\label{fig:geofwi-inits}
\end{figure}

\subsubsection{Out-of-distribution evaluation}
\label{sec:exp-geofwi-bp-ood}

The performance of all learning-based methods depends on how well the 
training set captures the features of the true and unknown medium. In Fig. 
\ref{fig:bp-ood-165} we illustrate the results for two $2.05\,$km$\times 3.0\,$km crops of the  challenging ``BP-2004 benchmark velocity model" \cite{billette2004}.
The crops contain heavy salt
structures (salt fraction $\approx 0.25$--$0.38$, vs $\approx 0.06$
for a representative GeoFWI salt model), placing them well outside
the GeoFWI training distribution. Not surprisingly, none of the learning-based methods do a good job of reconstructing $c(\bx)$. 

The ROM
approach gives reasonable results, but as explained in the next section, it carries a high computational cost.
This cost can be reduced by using the ROMNet output as an initial guess. That is to say, ROMNet 
alone may not be sufficient for giving accurate estimates of wave speeds that are far from those in the training set.
However, these estimates may be good starting values for the ROM-based inversion, as illustrated in Fig. \ref{fig:bp-ood-three-stage-evolution}.
We display there the estimated wave speed at $5, 10, 15$ and $20$ iterations of the Gauss-Newton method for the ROM-based inversion. 
The estimates are better than those in Fig. \ref{fig:bp-ood-165} (right column) obtained after $40$ iterations, because the starting guess is better. 

\begin{figure}[h]
\centering
\includegraphics[width=0.98\linewidth]{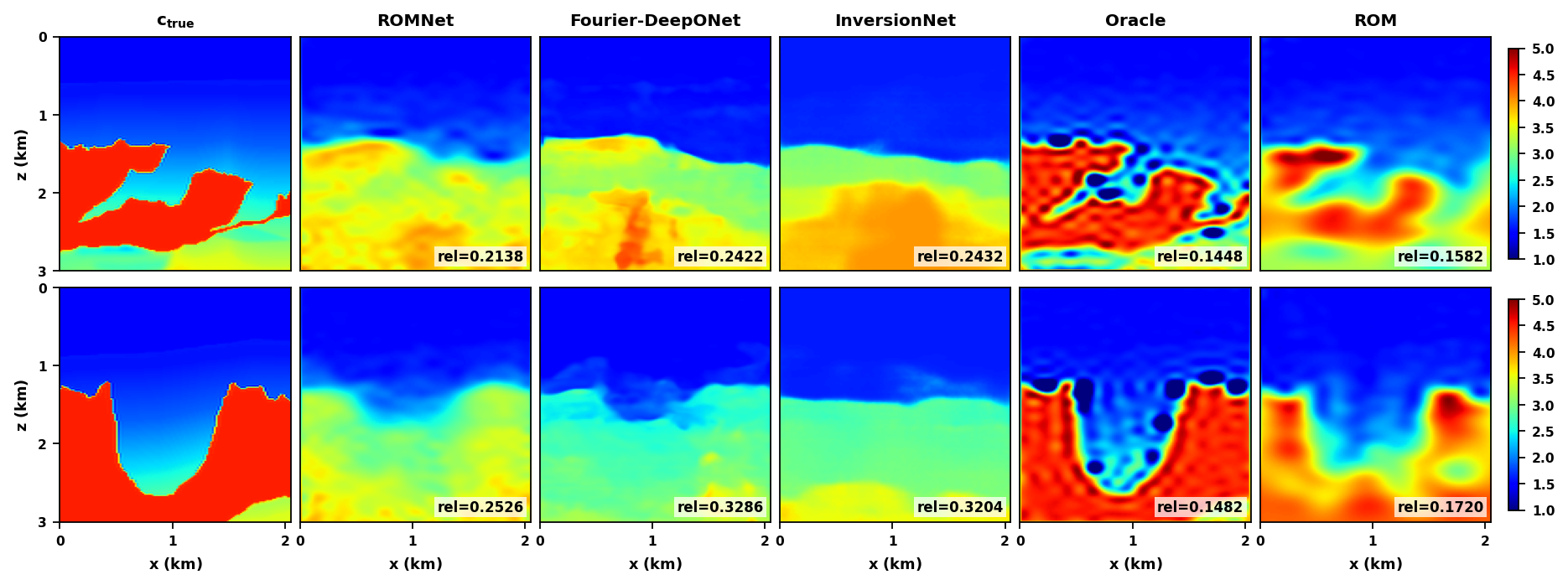}

\vspace{-0.1in}\caption{Two crops of the BP-2004 model. From left to right: True
$c(\bx)$, ROMNet, Fourier-DeepONet,
InversionNet, Oracle,
ROM.}
\label{fig:bp-ood-165}
\end{figure}

\begin{figure}[h]
\centering
\includegraphics[width=0.98\linewidth]{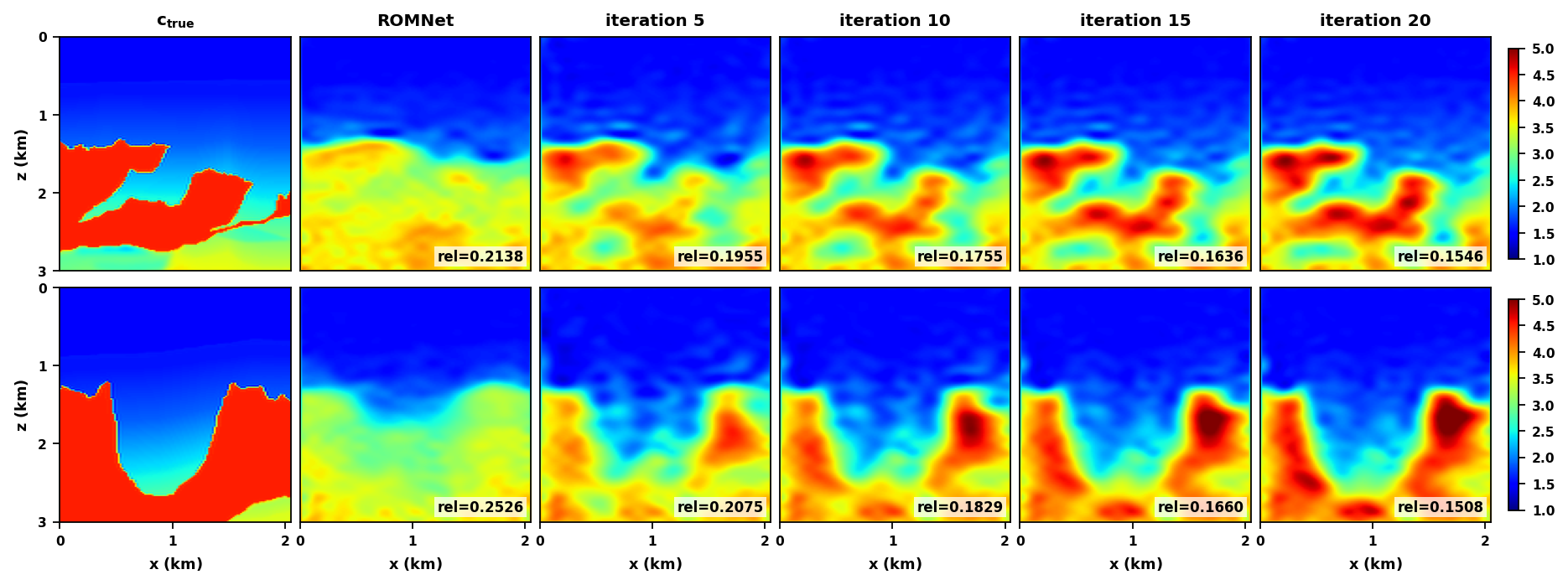}

\vspace{-0.1in}\caption{The evolution of the estimated wave speed with the ROM-based inversion, when starting with the initial guess given by ROMNet.
{From left to right: True $c(\bx)$, the ROMNet estimate (the same as in Fig. \ref{fig:bp-ood-165}) and the estimates after $5$, $10$, $15$ and $20$ Gauss-Newton iterations.
}
}
\label{fig:bp-ood-three-stage-evolution}
\end{figure}

\subsection{Computational cost and acceleration}
\label{sec:exp-speed}
Both the ROM inversion and ROMNet require computing from the data 
the operator ROM matrix $\bA^{\RM}$. The most expensive part of this calculation is the 
block-Cholesky factorization of $\bM$, which requires $O\big((n m)^3 \big)$ operations.
This cost can be controlled by time windowing the data (reducing $n$) and segmenting the array in sub-apertures
(reducing $m$). 

The optimization carried out with the Gauss-Newton iteration  is 
described in appendix \ref{sect:GaussN}. The cost difference between the ROM inversion and ROMNet lies in the 
computation of the Jacobian and the objective function. The
 ROM inversion must re-solve the wave equation for all $m$ sources at
every iteration to build $\bA^{\RM}(\hat c(\,\cdot\,;\bet))$. The cost  is $O(m\,n_{\rm fine}N_{\rm fine})$ per solve,
where $n_{\rm fine}$ and $N_{\rm fine}$ denote the number of time steps and points on the fine temporal and spatial grids
used by the solver. We used 
\[
n_{\rm fine} = (2n-1)20 + 1 \quad \mbox{and} \quad 
N_{\rm fine} = (N_y-2)(N_x-2),
\]
with $n$, $N_x$ and $N_y$ given in Tables \ref{tab:setup-rg} and \ref{tab:setup-geofwi} in appendix \ref{sec:exp-setup}.
%We timed the matrix-free variant of the computation of the Jacobian, see appendix \ref{sect:GaussN} for its description.

The ROMNet approach evaluates the Stage~1 network once. There is no need to solve 
the wave equation. Each iteration assembles the matrix 
$\bA^{\ML}(\hat c(\,\cdot\,;\bet))$ according to the formula \eqref{eq:3.2}. 
The cost per iteration is that of numerical quadrature i.e.,  $O\!\big(N_{\rm fine}(nm)^2\big)$.
The Jacobian is obtained by automatic differentiation and, as for the ROM inversion, it is applied matrix-free, with the regularized normal equations solved by conjugate gradients (see appendix \ref{sect:GaussN}).

Table~\ref{tab:speed} confirms,  on a single NVIDIA A100 (80\,GB) GPU,  that the dominant 
cost comes from the wave equation solver.  The ROM inversion costs $19.4$\,s per iteration while
ROMNet costs $0.26$\,s per iteration. End to end, for an inversion with $K=40$ iterations, the cost drops from $14$\,min
to $10.5$\,s. 
The one-time Stage~1 network evaluation adds a negligible time cost of $6$\,ms.

\begin{table}[h]
\centering\small
\begin{tabular}{lccc}
\toprule
 & per iter. & total ($K\!=\!40$) & error \\
\midrule
 ROM & $19.4$\,s & $14$\,min  & $0.0295$ \\
ROMNet & $\mathbf{0.26}$\,s & $\mathbf{10.5}$\,s    & $0.0296$ \\
\bottomrule
\end{tabular}
\vspace{0.02in}
\caption{Wall time and $L^2$ relative error for ROM and ROMNet. Per-iteration times exclude the one-time operator build and
compilation, and the totals are measured end-to-end sums, averaged over the three test samples on an idle GPU. {The errors are averaged over the same three media, which is why they differ slightly from those of Table~\ref{tab:indist-per-sample}, and the ROM inversion and ROMNet are both run for $K = 40$ iterations.} The two methods use the same configuration described in appendix \ref{sect:GaussN}.}
\label{tab:speed}
\end{table}

Fig.~\ref{fig:rg-speed-rel-convergence-p785} tracks the $L^2$ relative error
of the estimated wave speed at every iteration (mean over three examples), starting with the same initial guess $c_o$. ROMNet reaches its
plateau at about $15$ iterations. The ROM error decreases monotonically to a comparable level at iteration $40$, at about $75$ times the per-iteration cost of ROMNet.

\begin{figure}[t]
\centering
\includegraphics[width=0.45\linewidth]{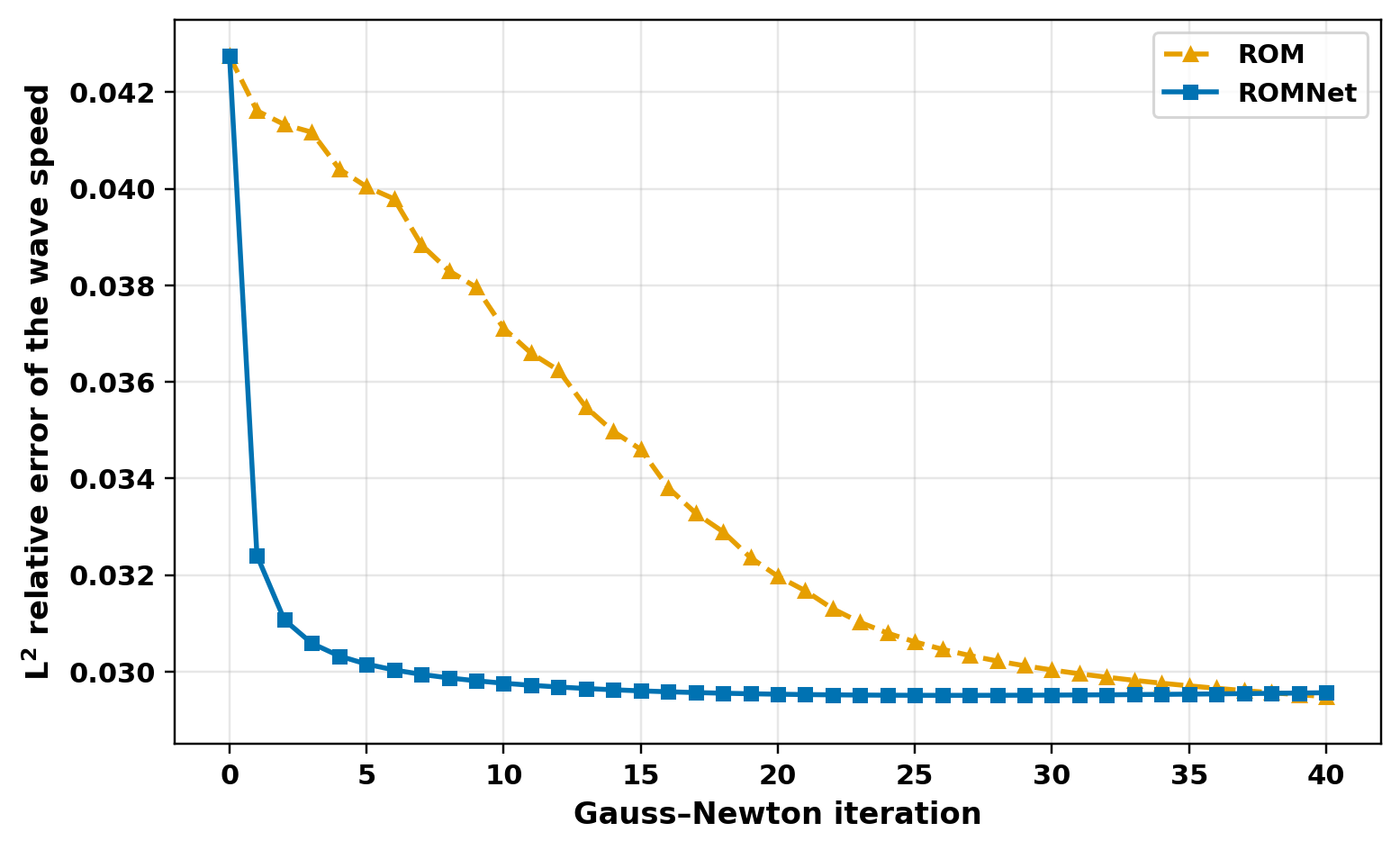}

\vspace{-0.1in}\caption{Random
Gaussians dataset: $L^2$ relative error of the reconstruction per Gauss--Newton iteration, for the same configuration as in Table~\ref{tab:speed}.}
\label{fig:rg-speed-rel-convergence-p785}
\end{figure}

\subsection{Dealing with noise}
\label{sect:noise}
If the data are noisy, the mass matrix $\bM$, with blocks computed as in equation \eqref{eq:2.18}, will not be symmetric and positive definite.
Symmetry is easy to fix: Either we take $(\bM + \bM^T)/2$ or, we compute the blocks of $\bM$ 
in the block upper triangular part and then set the blocks in the lower part by symmetry. The positive definiteness must be enforced with careful regularization, because the success of ROMNet relies on having the correct algebraic structure of $\bA^{\RM}$. One regularization 
procedure that gives such a structure was introduced in \cite[Appendix E]{Borcea2022rom}. It involves projections and orthogonal transformations that are medium dependent, and are thus difficult to use  in the training of the 
neural network. Here we use a simpler regularization,  based on the following observation:  All the diagonal blocks of $\bM$ contain the 
term $\bcD(0)/2$ which is computed (not measured) in the medium with constant speed $\bar{c}$. By ``boosting" this noiseless data matrix, as in $(1+\epsilon)\bcD(0)$, with user defined $\epsilon>0$ depending on the noise level, we are adding to $\bM$ the positive definite\footnote{The matrix $\bcD(0)$ is positive definite according to definition \eqref{eq:2.4}, since it is the Gramian of the $m$ snapshots gathered in $\bu_0(\bx)$.}  block diagonal matrix {$\frac{\epsilon}{2} \mbox{diag} \Big(2\bcD(0), \bcD(0), \ldots, \bcD(0)\Big)$, while maintaining the desired Toeplitz+Hankel block algebraic structure of the mass matrix. The addition shifts the spectrum of
\begin{equation}
\bM^\epsilon = \bM + \frac{\epsilon}{2} \mbox{diag} \Big({2\bcD(0), }\bcD(0), \ldots, \bcD(0)\Big)
\label{eq:regM}
\end{equation}
to positive eigenvalues, so we can compute its block Cholesky square root  and then its inverse.

The regularized ROMNet procedure remains as in Algorithm \ref{alg:two-stage}, except that at line~\ref{ln:computeROM} we compute the 
ROM from the regularized mass matrix \eqref{eq:regM} and the stiffness matrix modified accordingly to 
\begin{equation}
\bS^\epsilon = \bS - \frac{\epsilon}{2} \mbox{diag} \Big({2\ddot\bcD(0), }\ddot\bcD(0), \ldots, \ddot\bcD(0)\Big).
\label{eq:regS}
\end{equation}
Here we recalled the data driven calculation of $\bS$ given in equation \eqref{eq:2.20}. 

The training of the neural network is  modified by using the same zero-time boosting of the data. The inputs $\bA^{\RM}(c_i)$ are computed from simulated data contaminated with noise, at the same level as the measurements. In Fig. \ref{fig:noise-geofwi}  we used $\epsilon = 0.1$ and 
\begin{equation}
\bcD^{\rm noisy}(t_k;c_i) = \bcD(t_k;c_i) + \big[{\boldsymbol{\mathcal Z}}(t_k) + {\boldsymbol{\mathcal Z}}^T(t_k)\big]/2, \quad \mbox{for} ~k > 0, 
\label{eq:noise}
\end{equation}
where the entries of the $m\times m$ matrices ${\boldsymbol{\mathcal Z}}(t_k)$ are independent, identically distributed Gaussian random variables with mean zero and standard deviation $0.01 \, \beta$. Here  $\beta^2 = \frac{1}{m^2 N_t} \sum_{k=1}^{N_t}\|\bcD(t_k)\|_F^2$ is the mean square of the noiseless data entries, so the noise level is $1\%$. 
The noisy data are processed exactly as in \cite[Appendix A]{Borcea2022rom}.

Fig. \ref{fig:noise-geofwi} shows that for the three test media, the regularized ROMNet estimates from noisy data have practically the same quality as those from noiseless data.

\begin{figure}[t]
\centering
\includegraphics[width=0.6\linewidth]{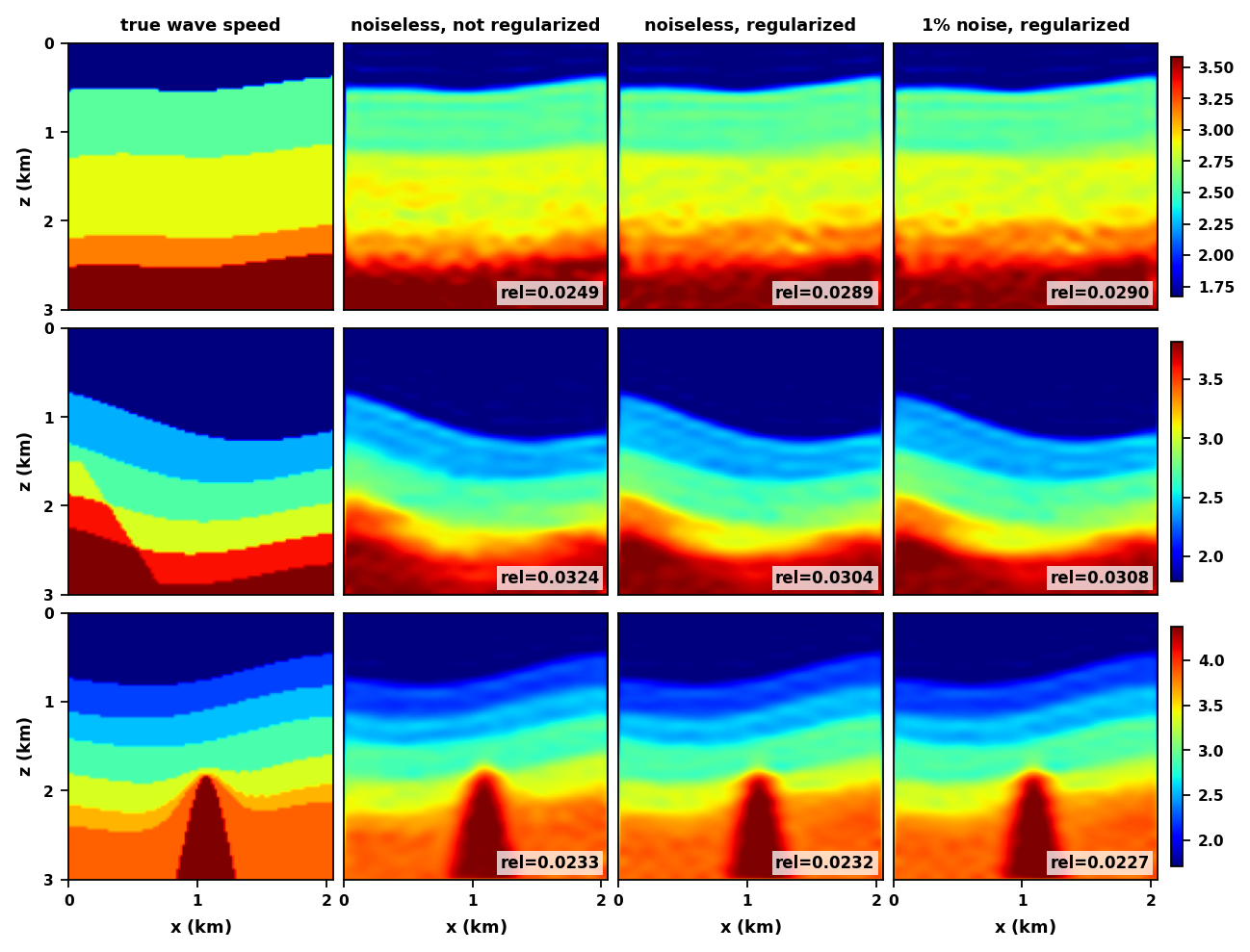}

\vspace{-0.1in} \caption{{GeoFWI dataset: ROMNet estimates for 3 test media, from noiseless data without and with the regularization \eqref{eq:regM}--\eqref{eq:regS}, and from data with $1\%$ additive noise, with the regularization.}}
\label{fig:noise-geofwi}
\end{figure}

\section{Summary}
\label{sect:4}
We introduced a novel approach to waveform inversion,  which seeks to estimate the unknown wave speed $c(\bx)$ of a medium from backscattering data gathered by an active array of sensors. It is a hybrid 
approach that combines recent advances in waveform inversion via data driven reduced order modeling and machine learning, hence the name ROMNet. The reduced order model (ROM) is an algebraic (matrix) surrogate of the wave equation operator. It is related 
to time snapshots of the data by a  nonlinear mapping that is well understood and is computed  in a mathematically tractable, non-iterative way. The estimation of $c(\bx)$ from the ROM is not well understood and so far, it has been formulated as a computationally expensive, nonlinear optimization. The main idea in this paper is to use a neural network to map the ROM to another matrix that has a simpler and explicit (quadratic) 
dependence on $c(\bx)$. The estimation of $c(\bx)$ from the network prediction is still formulated as an optimization problem. However,
this problem is easier to solve, and its computational cost is orders of magnitude less than that of the ROM-based waveform inversion.

We formulate and motivate the ROMNet approach starting from the wave equation and assess its performance with numerical simulations. The assessment includes a comparison with two commonly used learning-based data fitting approaches to waveform inversion: ``Fourier-DeepONet" and ``InversionNet". These approaches are trained on the same data sets as ROMNet. We used two such sets: The first consists of random media with variations of the wave speed modeled by a superposition of Gaussian functions with random amplitudes and standard deviations. The second data set is the publicly available GeoFWI data set. The result shows that ROMNet is the most accurate of the learning-based methods. The in-training-distribution results are consistently better and some out-of-distribution results, on known challenging problems, are good as well, slightly worse than those of the expensive, pure ROM-based inversion. The results also show that ROMNet,   like all other learning-based methods, cannot accurately estimate wave speeds far from those used during training. However, the estimates are good enough to provide otherwise unavailable good initial guesses for the pure ROM-based inversion. Such guesses lead to convergence in fewer iterations and, therefore, a significant reduction in computational cost.

\section*{Acknowledgment}
This work relates to the Air Force award FA9550261B167, issued by the Air
Force Office of Scientific Research, to the Department of Navy award N000142612023, issued by the Office of
Naval Research, and to NSF CAIG RISE2530596.

\appendix
\section{Details of the numerical implementation}
\label{ap:A}
Here we give the details of our simulations. 

\subsection{Setup}
\label{sec:exp-setup}
The wave equation solver is as in \cite{Borcea2022rom}. The 
Laplacian is discretized with a five point finite difference stencil on a rectangular grid described in Tables \ref{tab:setup-rg} and \ref{tab:setup-geofwi}. 
The domain $\Omega$ is rectangular, and we use homogeneous Dirichlet boundary conditions. The second time derivative is 
approximated by a three-point finite difference scheme, on a fine time grid with step 
$\tau/20$. {The wave equation solver, the computation of the operator ROM, and the Gauss-Newton inversions, including the Jacobians obtained by automatic differentiation, are implemented in JAX. The neural network and its training are implemented in PyTorch.}

The array response matrix is given by the numerically approximated 
pressure, evaluated at the sensor locations, as given in equation \eqref{eq:I3}.
The data matrices are computed as in \eqref{eq:2.3}, where the convolution is evaluated using the 
Fast Fourier transform. The array is linear and placed horizontally, near the top boundary of $\Omega$. The sensors  are uniformly separated in the aperture $a$ of the array.  Both $a$ and the distance (depth) of the array from $\partial \Omega$ are indicated in Tables~\ref{tab:setup-rg} and \ref{tab:setup-geofwi}. 

The search dictionary  $\bPhi(\bx)$  in the parametrization \eqref{eq:searchC} consists of Gaussians $\phi_l(\bx) =\exp\big[-(x - x_l)^2/(2\rho_x^2) - (y-y_l)^2/(2\rho_y^2)\big]$, with centers $(x_l,y_l)$ on a uniform grid covering $\Omega$. See Table \ref{tab:setup-inversion} for the number $N_c$ of functions and the grid points.

\begin{table}[h]
\centering
\small\color{black}
\begin{tabular}{lp{9.3cm}}
\toprule
\multicolumn{2}{c}{\textbf{Random Gaussians}} \\
\midrule
Domain, grid & $\Omega = [0,2.0]\times[0,2.5]$km, uniform grid of $100 \times 125$ points \\
Probing pulse & $\om_o/(2\pi) = 6$Hz, $B/(2\pi) = 4$Hz \\
$\bar{c}$ and $\la_{\rm min}$ & $\bar{c} = 3$km/s, $\la_{\rm min} = 2\pi \bar{c}/(\om_o + B) = 0.3$km \\
Array & $m = 10$ sensors at depth $0.14$km, at horizontal positions $0.30, 0.456, \ldots, 1.70$km (spacing $0.156$km, aperture $a = 1.4$km) \\
Time sampling & $\tau = 0.040$s, $n = 16$ snapshots \\
Reference speed &$c_o(\bx) \equiv \bar{c} = 3$km/s \\
Wave speed models & $c(\bx) = \mathrm{clip}_{[1.6,4.0]}\Big(\bar{c} + \sum_{p} \eta_p\, \kappa_{\rho_p}(\bx - \bx_p)\Big)$, {with the sum over the grid points $\bx_p$; $\kappa_\rho(\bx) \propto e^{-|\bx|^2/(2 \rho^2)}$, normalized to unit sum over the grid.} The weights $\eta_p$ are drawn uniformly and independently from $[-4,4]$km/s and multiplied by a taper that vanishes within $0.375$km of the top and bottom boundaries and $0.456$km of the lateral boundaries (buffer zone, containing the array). $\rho_p$ drawn independently, with equal probability, from $5$ values log-spaced in $[\la_{\rm min}/9, \la_{\rm min}] = [0.033, 0.3]$km. $\mathrm{Clip}$ is cutoff to $[1.6,4.0]$km/s \\
Data sets & $10000$ media: $7000$ / $1500$ / $1500$ (training / validation / test) \\
\bottomrule
\end{tabular}
\caption{Setup for the Random Gaussians data set.}
\label{tab:setup-rg}
\end{table}

\begin{table}[h]
\centering
\small\color{black}
\begin{tabular}{lp{9.3cm}}
\toprule
\multicolumn{2}{c}{\textbf{GeoFWI \cite{li2026geofwi}}} \\
\midrule
Domain, grid & $\Omega = [0,2.05]\times[0,3.0]$km, uniform grid of $123 \times 180$ points. GeoFWI models, given on $100 \times 100$ grids, are resampled to this grid. \\
Probing pulse & $\om_o/(2\pi) = 6$Hz, $B/(2\pi) = 4$Hz \\
$\bar{c}$ and $\la_{\rm min}$ & $\bar{c} = 1.5$km/s (the smallest wave speed in the data set, in water), $\la_{\rm min} = 2\pi \bar{c}/(\om_o + B) = 0.15$km \\
Array & $m = 10$ sensors at depth $0.15$km, at horizontal positions $0.30, 0.461, \ldots, 1.75$km (spacing $0.161$km, aperture $a = 1.45$km) \\
Time sampling & $\tau = 0.0435$s, $n = 40$ snapshots for ROMNet; the ROM inversion uses $n = 56$ snapshots, with the domain extended below the imaging region, so that the bottom boundary has no effect. \\
Reference $c_o(\bx)$ & the per-depth median of the training wave speeds; it equals $1.69$km/s at the depth of the array and $3.83$km/s at the bottom \\
Wave speed models & GeoFWI media with folds, faults and salt bodies, clipped to $[1.5, 4.5]$km/s; the water-shallow selection \eqref{eq:water-shallow} keeps the media with wave speed below $2.0$km/s above the array \\
Data sets & water-shallow subset (\ref{eq:water-shallow}) split in $36000$ / $4000$ / $9476$ (training / validation / test)  \\
\bottomrule
\end{tabular}
\caption{ Setup for the GeoFWI data set.}
\label{tab:setup-geofwi}
\end{table}

\begin{table}[h]
\centering
\small\color{black}
\begin{tabular}{lp{4.2cm}p{4.2cm}}
\toprule
 & \textbf{Random Gaussians} & \textbf{GeoFWI} \\
\midrule
Number of functions & $28 \times 28 + 1 = N_c = 785$ & $20 \times 32 + 1 = N_c = 641$ \\
Spacing of the centers & $0.072$km $\times$ $0.090$km & $0.105$km $\times$ $0.095$km \\
$\rho_x, \rho_y$ ($0.7 \times$ spacing) & $0.050$km, $0.063$km & $0.074$km, $0.067$km \\
Windowing $\gamma$ & $5$ & $5$ \\
\bottomrule
\end{tabular}
\caption{Setup of the Stage 2 (Gauss--Newton) inversion for the two data sets.}
\label{tab:setup-inversion}
\end{table}

\subsection{Gauss-Newton iteration}
\label{sect:GaussN}
Following \cite{Borcea2022rom}, we use a Gauss-Newton iteration to solve the minimization problem \eqref{eq:gammaROM} with 
Tikhonov regularization i.e., the objective function is ${\mathcal{O}}^{\RM}_{\gamma}(\bet) + \mu \|\bet\|_2^2$. 
The regularization parameter $\mu$ is chosen as follows: Let 
\begin{equation}
\label{eq:GN1}
\boldsymbol{\mathcal{E}}(\bet^{(i)}) = W_{\gamma}\!\bigl[\bA^{\RM}\!\bigl(\hat c(\cdot; \bet^{(i)}) \bigr)
- \bA^{\RM}\bigr] \in \RR^{N_\gamma},
\end{equation}
be the residual vector, with $N_\gamma$ defined in equation \eqref{eq:Ngamma}. The Jacobian matrix is
\begin{equation}
\label{eq:GN2}
\boldsymbol{J}^{(i)} = \nabla_{\bet}\boldsymbol{\mathcal{E}}(\bet^{(i)}) \in \RR^{N_\gamma \times N_c},
\end{equation}
where typically, $N_c \ll N_{\gamma}$. {We take the same $\mu = 2 \times 10^{-3} \sigma_{\max}^2$ at all the iterations,  where $\sigma_{\max}$ is the largest singular value of the Jacobian $\boldsymbol{J}^{(0)}$ of the first iteration, estimated by eight power iterations.} \{Each Gauss-Newton step is adjusted using the same backtracking line search as in Algorithm \ref{alg:two-stage}: the step length starts at $\alpha = 1$ and is halved until the objective function decreases.} 

{The Jacobian \eqref{eq:GN2} of the ROM inversion is never formed explicitly.} We apply $\bJ$ and $\bJ^\top$ through {jvp
and vjp}  (the Jacobian-vector product $\bJ \boldsymbol{v}$ and the vector-Jacobian product $\boldsymbol{w}^T \bJ$), computed by forward and reverse mode automatic differentiation of the map $\bet \mapsto \boldsymbol{\mathcal{E}}(\bet)$, inside a conjugate-gradient (CG) solve of the normal equations. This
is the cost-optimal form matching the {$O(mn_{\rm fine} N_{\rm fine} )$}
estimate of \cite[\S5]{Borcea2022rom}. Each CG step costs one jvp and one vjp, that is, two
wave-equation solves.

{The Gauss-Newton iteration for the ROMNet problem \eqref{eq:Odkhat} is similar to the above. The residual vector is $W_\gamma[\bA^{\ML}(\hat c(\cdot;\bet))] - \cG_{\btheta}\big[\bA^{\RM}\big]$, whose second term is fixed. No wave equation is solved: at each iteration, $\bA^{\ML}(\hat c(\cdot;\bet))$ is assembled by numerical quadrature from the definition \eqref{eq:3.2}, using $\bV(\bx;c_o)$ computed in the reference medium, which does not change along the iterations. The Jacobian is obtained by automatic differentiation of this quadrature and is applied in the same matrix-free way, with the line search described above and the same damping $\mu$.} 

The causal construction of the operator ROM and of the learned matrix $\bA^{\ML}$ is exploited in the Gauss-Newton iterations using time-windowing of the data. This allows  the estimation of $c(\bx)$ in a layer-peeling fashion, starting from vicinity of the array and moving deeper inside the medium. The time-windowing is implemented as in \cite[Algorithm~2]{Borcea2022rom}:
For the simulations trained on the random Gaussian data set, the time window grows from $2\tau$ to $16 \tau$, in steps of $2 \tau$. 
For each intermediary time window we perform 3 Gauss-Newton iterations. For the full time window we perform the remaining iterations, 
up to the total $K = 40$.
For the simulations trained on the GeoFWI data set, the time window grows from $14\tau$ to $56 \tau$, in steps of $7 \tau$. 
For each intermediary time window we perform 3 Gauss-Newton iterations. For the full time window we perform the remaining iterations, 
up to the total $K = 40$.

\subsection{Neural network and training}
\label{ap:A.L}
Here we describe briefly the architecture of the neural network and the learning procedure. Recall that the network and its training are implemented in PyTorch.

{\textbf{Data normalization:} The training set consists of the pairs $\big(\bA^{\RM}(c_i), \bA^{\ML}(c_i)\big)$, $i = 1, \ldots, N_{\rm train}$, computed for the media $c_i(\bx)$ of the training split (Tables \ref{tab:setup-rg} and \ref{tab:setup-geofwi}). The network works with the residual $\bA^{\RM}(c_i) - \bA^{\RM}(c_o)$ of the ROM with respect to the ROM of the reference medium, and with the windowed target $W_\gamma[\bA^{\ML}(c_i)]$. On the input side the network reads the whole residual matrix: its entries are collected in a vector of length $N_n = nm(nm+1)/2$ by the operator $W_n$.
This is the operator $W_\gamma$ of section \ref{sect:MLAlg}, with $\gamma = n$. It extracts the whole upper triangular part, which determines the matrix because it is symmetric. On the output side, only the entries in the first $\gamma m$ diagonals are kept, in a vector of length $N_\gamma$ extracted by $W_\gamma$. From the training set we compute once the entry-wise means $\mu_X$, $\mu_Y$ and standard deviations $\sigma_X$, $\sigma_Y$ of the vectors $W_n[\bA^{\RM}(c_i) - \bA^{\RM}(c_o)]$ and $W_\gamma[\bA^{\ML}(c_i)]$, $i = 1, \ldots, N_{\rm train}$, and store these four vectors. 

The normalized training data are 
\[
X_i = \big(W_n[\bA^{\RM}(c_i) - \bA^{\RM}(c_o)] - \mu_X\big) \oslash \sigma_X \in \RR^{N_n}, \]
and 
\[Y_i = \big(W_\gamma[\bA^{\ML}(c_i)] - \mu_Y\big) \oslash \sigma_Y \in \RR^{N_\gamma}, i = 1, \ldots, N_{\rm train},
\]
 where $\oslash$ denotes entry-wise division,  $X_i$ is the input of the network and $Y_i$ its target.}

\textbf{Network architecture:} The network $w_{\btheta}: \RR^{N_n} \to \RR^{N_\gamma}$ maps the normalized input $x$ to an approximation of the normalized target $y$. 

The input vector $X$ can be organized in blocks of rows. The $j^{\rm th}$ block corresponds to the entries in the $(j,j), (j,j+1), \ldots, (j,n-1)$ blocks of the residual. We pad it with zero blocks, so that all the block rows have the same length $n$. This way we obtain an array called token $j$. Therefore, the network treats $X$ as the sequence of $n$ tokens, ordered in time. Each token is embedded in a vector $h_j \in \RR^{1536}$ by a small network, with the same weights for all the tokens. The sequence $(h_0, \ldots, h_{n-1})$ is then processed by a causal temporal convolutional network that is described in \cite[Section~3]{bai2018empirical}. We use $8$ residual blocks with convolutions of kernel size $3$ along the token index, arranged so that the output $h_j'$ for token $j$ depends only on $h_i$ for $i \le j$. No dilation is used, because the sequence has only $n$ tokens. This respects the causality of the ROM: the first block rows of $\bA^{\RM}(c_i)$ depend on the medium down to the depth swept by the wave up to the corresponding time, and so do the matching block rows of $\bA^{\ML}(c_i)$ \cite{Borcea2022rom}. {Finally, another small network, with the same weights for all the tokens, maps $h_j' \in \RR^{1536}$ to the $\gamma$ blocks $(j,j), \ldots, (j,j+\gamma-1)$ of block row $j$ of $Y$
%, i.e., to the entries of block row $j$ in the first $\gamma m$ diagonals, 
and these are multiplied by the gate $g_j = \tanh\big(\|(\sigma_X \odot X + \mu_X)_j\|_{\rm rms}/\tau_g\big) \in [0,1)$, computed from the de-normalized input $\sigma_X \odot X + \mu_X$ (i.e., the residual $\bA^{\RM}(c_i) - \bA^{\RM}(c_o)$), where $(\cdot)_j$ denotes the entries in block row $j$, $\|\cdot\|_{\rm rms}$ is the Euclidean norm divided by the square root of the number of entries, $\odot$ denotes entry-wise multiplication, and $\tau_g = 0.1$. All the linear maps in the network have no bias, so the gate and the network output vanish when block row $j$ of the input ROM coincides with the reference ROM, and at the reference medium the prediction reduces to the mean $\mu_Y$ of the training targets.}

\textbf{Loss and optimization:} The loss is the mean squared error between $w_{\btheta}(X_i)$ and the target $Y_i$, over the $N_\gamma$ entries of the output band and the training set. In terms of the physical scale of the target it reads
\[
\mathcal{L}(\btheta) = \frac{1}{N_{\rm train}} \sum_{i=1}^{N_{\rm train}} \frac{1}{s_i^2} \big\| \sigma_Y \odot \big(w_{\btheta}(X_i) - Y_i\big) \big\|_2^2 ,
\]
where $\sigma_Y \odot (w_{\btheta}(X_i) - Y_i)$ is the error of the prediction of $W_\gamma[\bA^{\ML}(c_i)]$ and $s_i$ is the root mean square of the entries of the target residual $W_\gamma[\bA^{\ML}(c_i) - \bA^{\ML}(c_o)]$, so that all the media contribute comparably. A small quadratic penalty with weight $10^{-5}$ on the output of the network is added. The optimizer is AdamW \cite{loshchilov2017decoupled} with learning rate $2 \times 10^{-4}$, a cosine schedule, batch size $32$ and gradient clipping at unit norm; the training runs for $400$ epochs and the parameters with the smallest error on the validation set are kept.

{\textbf{Prediction:} The map $\cG_{\btheta}$ of \eqref{eq:NetLearn} is $w_{\btheta}$ composed with the normalization of the input and the de-normalization of the output, with the stored statistics: for a medium with operator ROM $\bA^{\RM}$, computed from the data,
\[
\cG_{\btheta}\big[\bA^{\RM}\big] = \mu_Y + \sigma_Y \odot w_{\btheta}(X), \qquad X = \big(W_n[\bA^{\RM} - \bA^{\RM}(c_o)] - \mu_X\big) \oslash \sigma_X .
\]
Here $\cG_{\btheta}\big[\bA^{\RM}\big]$ is the Stage 1 prediction of $W_\gamma[\bA^{\ML}]$ in Algorithm \ref{alg:two-stage}, used in \eqref{eq:Odkhat}.}

% --------------
\bibliographystyle{siam}
\bibliography{biblio}
\end{document}